\documentclass[reqno,11pt]{amsart}
\usepackage[letterpaper,hmargin=1.25in,vmargin=0.9in]{geometry}
\usepackage[T1]{fontenc}
\usepackage{lmodern}
\usepackage{microtype}
\usepackage{mathtools}
\usepackage{amssymb}

\usepackage{xcolor}
\usepackage{booktabs,array}
\usepackage{tikz}
\usepackage[unicode]{hyperref}
\hypersetup{
  colorlinks=false,
  citecolor=blue,
  linkcolor=blue,
  urlcolor=blue
}

\allowdisplaybreaks[1]
\renewcommand{\footnoterule}{\kern-3pt\hrule width5pc height0.4pt\kern2.6pt}

\newcommand{\R}{\mathbb{R}}
\newcommand{\Z}{\mathbb{Z}}

\newcommand{\doi}[1]{DOI: \href{https://doi.org/#1}{#1}}

\DeclareMathOperator{\conv}{conv}

\DeclareMathOperator{\vol}{vol}

\theoremstyle{plain}
\newtheorem{theorem}{Theorem}[section]
\newtheorem{proposition}[theorem]{Proposition}
\newtheorem{lemma}[theorem]{Lemma}
\newtheorem{corollary}[theorem]{Corollary}

\numberwithin{equation}{section}

\title[Lattice coverings by octahedra]{The lattice covering density of the regular octahedron}

\author[Yanlu Lian]{Yanlu Lian\textsuperscript{1,2}}
\address{\textsuperscript{1}School of Mathematics, Hangzhou Normal University, Hangzhou 311121, China}
\address{\textsuperscript{2}School of Computational and Mathematical Sciences, Cardiff University, Cardiff, UK}
\email{yllian@hznu.edu.cn}

\author[Fei Xue]{Fei Xue\textsuperscript{3,*}}
\thanks{\textsuperscript{*}Corresponding author.}
\address{\textsuperscript{3}School of Mathematical Sciences, Nanjing Normal University, Nanjing 210046, China}
\email{05429@njnu.edu.cn}

\thanks{This work is supported by the National Natural Science Foundation of China grant (NSFC12671118), the Natural Science Foundation of Zhejiang Province grant (LQ24A010008) and the Scientific Research Fund of Zhejiang Provincial Education Department grant (Y202351675).}

\hypersetup{
  pdftitle={The lattice covering density of the regular octahedron},
  pdfauthor={Yanlu Lian and Fei Xue},
  pdfsubject={Mathematics},
  pdfkeywords={Lattice covering density, regular octahedron, lattice tetrahedra, Kuhn triangulation, Abelian Cayley graphs}
}

\begin{document}

\begin{abstract}
We prove that the lattice covering density of the regular octahedron is $9/8$, establishing the optimality of the classical body-centered cubic covering. A hypothetical lattice covering of smaller density would yield a lattice basis whose six Kuhn tetrahedra each fit in a translate of the octahedron. A sharp determinant bound then gives a contradiction. As a consequence, the maximum order of an undirected Abelian Cayley graph with three generators and their inverses and diameter at most $k$ is $\frac{32}{27}k^3+O(k^2)$ as $k\to\infty$.
\end{abstract}

\maketitle

\section{Introduction}

Determining the lattice covering density of a convex body is a classical problem in the geometry of numbers and discrete geometry. For a convex body $C\subset\R^3$, its lattice covering density is
$$\theta^l(C)=\inf\left\{\frac{\vol(C)}{\det L}: L\text{ is a full-rank lattice and }C+L=\R^3\right\}.$$
Here $\det L$ is the volume of a fundamental parallelepiped, and $\vol(C)/\det L$ is the average covering multiplicity. Even in three dimensions, determining this infimum requires control of all covering lattices, not only perturbations of a candidate optimum. For the Euclidean unit ball $B^3$, Bambah~\cite{Bambah1954} proved in 1954 that
$$\theta^l(B^3)=\frac{5\sqrt5\pi}{24}\approx 1.4635,$$
attained by a body-centered cubic lattice; see~\cite{Barnes1956,Few1956}.

We study the corresponding problem for the regular octahedron
$$O=\{\mathbf{x}\in\R^3:|x_1|+|x_2|+|x_3|\le 1 \},$$
the unit ball of the $\ell_1$ norm. A classical body-centered cubic construction gives a lattice covering of $O$ with density $9/8$. Geometrically, the construction can be obtained from a tiling by truncated octahedra. The lattice-tiling property of this parallelohedron goes back to Fedorov's work of 1885; see the history in~\cite{LLFZ2024}. We give a direct verification of the covering in Section~\ref{sec:proof-outline}. Dougherty and Faber studied this construction through the degree-diameter problem for Abelian Cayley graphs. They proved its strict local optimality and conjectured its global optimality, equivalently that the optimal lattice covering efficiency is $8/9$; see~\cite[Theorem~7.2 and the preceding discussion]{DF2004}. Their local result does not exclude a more efficient covering by a lattice far from the body-centered cubic lattice. They also describe a reduction of the global problem to finitely many optimization problems, each maximizing a cubic determinant subject to linear inequalities. The difficulty is the large number of possible covering configurations~\cite{DF2004}.

Let $\theta^t(O)$ denote the translative covering density of $O$, where the centers need not form a lattice. Since every lattice covering is a translative covering, the lower bound of Li, Lian, Fu, and Zhang~\cite{LLFZ2024}, together with the classical body-centered cubic construction, gives
$$1+6.6\times10^{-8}\leq \theta^t(O) \leq \theta^l(O)\leq \frac{9}{8}.$$
Li and Zong~\cite{LZ2025} determined the simultaneous lattice packing-covering constant of $O$ to be $7/6$. This problem additionally requires the same lattice to give a packing by $O$ and a covering by a dilate of $O$.

We prove that the classical body-centered cubic covering is globally optimal among lattice coverings of the regular octahedron.
\begin{theorem}\label{thm:main}
The lattice covering density of the regular octahedron is
$$\theta^l(O)=\frac{9}{8}.$$
\end{theorem}

A known sufficient condition for $C+L=\R^3$ is that $C$ contain a translate of each of the six Kuhn tetrahedra associated with one basis of $L$~\cite{Balitskiy2020,LX2026}. We show that a hypothetical lattice covering of $O$ with density less than $9/8$ would yield a basis of its covering lattice for which these six containment conditions hold. A sharp determinant bound for this basis then gives a contradiction. Section~\ref{sec:proof-outline} explains this reduction.

The theorem also determines the asymptotic maximum order of Abelian Cayley graphs with three generators and their inverses as the diameter grows. Let $M_3(k)$ denote this maximum order when the diameter is at most $k$. Walks of length at most $k$ are represented by integer vectors in $kO$. Through the relation between graph diameters and lattice coverings established by Dougherty and Faber~\cite{DF2004}, the identity $1/\theta^l(O)=8/9$ gives
$$M_3(k)=\frac{32}{27}k^3+O(k^2) \qquad(k\to\infty).$$
See Corollary~\ref{cor:graph-asymptotic}.

\section{Normalization and proof outline}\label{sec:proof-outline}

\subsection{The determinant problem}

For $\mathbf{x}=(x_1,x_2,x_3)\in\R^3$, write $|\mathbf{x}|_1=|x_1|+|x_2|+|x_3|$ and set
$$K=\frac{3}{2}O,\qquad \vol(K)=\frac{9}{2}.$$
The lower bound in Theorem~\ref{thm:main} is equivalent to
\begin{equation}\label{eq:determinant-goal}
  K+L=\R^3\Longrightarrow\det L\leq 4.
\end{equation}
The body-centered cubic lattice
\begin{equation}\label{eq:bcc-lattice}
  L_0=2\Z^3\cup\bigl((1,1,1)+2\Z^3\bigr)
\end{equation}
has determinant $4$, since $[L_0:2\Z^3]=2$. For $\mathbf{x}\in\R^3$, let $\delta_i\in[0,1]$ be the distance from $x_i$ to the nearest even integer. Its distance to the nearest odd integer is $1-\delta_i$, so
$$\operatorname{dist}_1(\mathbf{x},L_0) =\min\left\{\sum_{i=1}^3\delta_i,\;3-\sum_{i=1}^3\delta_i\right\} \leq\frac{3}{2}.$$
Thus $K+L_0=\R^3$, with density $(9/2)/4=9/8$.

\smallskip\noindent
{\it Notation.} We use $[\mathbf{a}~\mathbf{b}~\mathbf{c}]$ for the matrix with these columns and $[\mathbf{a};\mathbf{b};\mathbf{c}]$ for the matrix with these rows. Volumes are three-dimensional unless another dimension is specified. The vectors $\mathbf{e}_1,\mathbf{e}_2,\mathbf{e}_3$ form the standard basis of $\R^3$, and $\mathfrak{S}_3$ denotes the permutations of three indices. Integer coordinates will always be taken relative to a basis matrix $B$, whose columns form a basis of $L$. Thus $\mathbf{v}\in\Z^3$ represents the lattice vector $B\mathbf{v}$. A change of these coordinates leaves the octahedron $K$ fixed.

\subsection{Outline of the proof}

The construction above gives the upper bound $9/8$. For the reverse inequality, suppose that
\begin{equation}\label{eq:counterexample-assumption}
  K+L=\R^3,\qquad D:=\det L>4.
\end{equation}
We will find a basis matrix $B$ of $L$ such that the six sets
\begin{equation}\label{eq:intro-six-paths}
  \Delta_{ijk} =\{\mathbf0,\mathbf e_i,\mathbf e_i+\mathbf e_j, \mathbf e_i+\mathbf e_j+\mathbf e_k\}, \qquad (i,j,k)\in\mathfrak S_3,
\end{equation}
satisfy
$$B\Delta_{ijk}\subseteq\mathbf t_{ijk}+K \qquad ((i,j,k)\in\mathfrak S_3).$$
The convex hulls of the six sets $\Delta_{ijk}$ form the standard Kuhn triangulation of $[0,1]^3$~\cite{Kuhn1960}. The six containing translates may differ. These containments are known to suffice for a lattice covering; see~\cite{Balitskiy2020} and~\cite[Theorem~1.1]{LX2026}. Here we derive them from the covering assumption under $D>4$.

Call lattice points {\it co-covering} if they lie in a common translate of $K$. Section~\ref{sec:lattice-topology} constructs a finite integral periodic chain $z$ of co-covering formal tetrahedra satisfying
$$\partial z=0,\qquad d(z)=6.$$
The first identity expresses cancellation of oriented triangular faces, with common lattice translations identified. The second expresses total signed volume $D$ and ensures that some tetrahedron is nondegenerate.

Every edge vector of a co-covering tetrahedron belongs to $K-K=3O$. Call a nonzero lattice vector $\mathbf v\in L$ {\it short} if $|\mathbf v|_1\le3$, and regard the pair $\{\mathbf v,-\mathbf v\}$ as one {\it direction}. The determinant estimate in Section~\ref{sec:four-vector} shows that every nondegenerate tetrahedron with short edges is {\it unimodular}, meaning that its three edge vectors from any vertex form a basis of $L$. Section~\ref{sec:short-vectors} proves that every nonzero lattice vector has $\ell_1$ norm greater than $3/2$ and that each nonzero class of $L/2L$ contains at most two short directions. The first restriction also implies that every short vector is {\it primitive}, meaning that it is not a nontrivial integer multiple of another lattice vector.

The determinant bound also controls a point joined by short vectors to the three vertices of a fixed triangular face. Lemma~\ref{lem:degenerate-bridge} gives explicit bounds on its coordinates in the basis supplied by a nondegenerate tetrahedron. If the point lies in the face plane, the resulting four-point set must be a parallelogram. Subtracting the reflected covering chain cancels these parallelogram terms. The remaining degenerate terms share no face with the nondegenerate terms, so the latter still have cancelling faces. Proposition~\ref{prop:closed-component} then selects a chain $c$ of nondegenerate co-covering tetrahedra with face-connected support and
$$\partial c=0,\qquad d(c)\neq 0.$$

For a point outside the face plane, the same coordinate bounds give the ten possible apices in Lemma~\ref{lem:shared-face-apices}. In particular, tetrahedra sharing a face have at least four common edge directions. Lemma~\ref{lem:two-tetrahedra} classifies the union of two different direction sets with this overlap. When the union has eight directions, the determinant and parity restrictions show that these same directions suffice for every tetrahedron reached through shared faces. Families with only six or seven directions are treated by a direct choice of basis. Thus Lemma~\ref{lem:connected-eight-directions} places all edge vectors of $c$ in one common set of eight directions.

With these coordinates fixed, Section~\ref{sec:manual-topology} lists the possible tetrahedron translation classes. Three face cancellation equations, together with $d(c)\ne0$, force a complete family of six Kuhn tetrahedra for one lattice basis. The tetrahedra obtained in this way are co-covering, so they give the six required containments.

Finally, Section~\ref{sec:six-path} converts the containment conditions into inequalities for a transformed basis matrix. Known descriptions of the Clifford and Mermin polytopes~\cite{BCLLSU2006,OCI2023} yield a bound on the sum of its singular values. The arithmetic--geometric mean inequality then gives
$$|\det B|\le\frac{32}{27}r^3$$
whenever the six sets $B\Delta_{ijk}$ lie in translates of $rO$. Taking $r=3/2$ gives
$$D=|\det B|\leq\frac{32}{27}\left(\frac{3}{2}\right)^3=4,$$
contrary to~\eqref{eq:counterexample-assumption}. Section~\ref{sec:verification} completes the proof and derives the consequence for Abelian Cayley graphs.

\section{A periodic chain of co-covering tetrahedra}\label{sec:lattice-topology}

Let $L$ satisfy $K+L=\R^3$, and write $D=\det L$. The goal of this section is Proposition~\ref{prop:nerve-degree}. It constructs a finite integral periodic chain of co-covering tetrahedra whose oriented faces cancel and whose total signed volume is $D$. We obtain this chain by replacing the vertices of a sufficiently fine periodic triangulation with covering lattice centers.

\subsection{Co-covering lattice points}

A finite set $F=\{\boldsymbol\ell_0,\ldots,\boldsymbol\ell_k\}\subset L$ is {\it co-covering} if the translates $K+\boldsymbol\ell_i$ have a common point. Since $K=-K$, this is equivalent to
$$\bigcap_{i=0}^k(K+\boldsymbol\ell_i)\ne\varnothing \quad\Longleftrightarrow\quad F\subset\mathbf{x}+K\quad\text{for some }\mathbf{x}\in\R^3.$$
We shall triangulate space into small tetrahedra and replace each vertex $\mathbf{x}$ by a lattice point $\boldsymbol\lambda$ with $\mathbf{x}\in K+\boldsymbol\lambda$. The next lemma ensures that the four chosen lattice points for each small tetrahedron are co-covering.

\begin{lemma}\label{lem:intersection-thickening}
For the fixed lattice $L$, there is $\varepsilon=\varepsilon(L)>0$ with the following property. Let $k\ge0$ and let $\mathbf{x}_0,\ldots,\mathbf{x}_k$ have pairwise $\ell_1$-distance less than $\varepsilon$. If $\boldsymbol\lambda_0,\ldots,\boldsymbol\lambda_k\in L$ satisfy $\mathbf{x}_i\in K+\boldsymbol\lambda_i$ for every $i$, then
$$\bigcap_{i=0}^k(K+\boldsymbol\lambda_i)\ne\varnothing.$$
\end{lemma}

\begin{proof}
Put $r=3/2$ and
$$E=\{\mathbf{v}\in L:|\mathbf{v}|_1\leq 2(r+1)\}.$$
For each $A\subseteq E$ containing $\mathbf{0}$, let
$$\rho(A)=\min_{\mathbf{x}\in\R^3} \max_{\boldsymbol\ell\in A}|\mathbf{x}-\boldsymbol\ell|_1.$$
This is the radius of the smallest $\ell_1$ ball containing $A$; thus $\rho(A)\le r$ exactly when $A$ is co-covering. The minimum exists because the maximum is continuous and tends to infinity with $|\mathbf{x}|_1$. The finite set $E$ has only finitely many subsets, so the radii greater than $r$ are separated from $r$ by a positive gap. We may therefore choose
$$\varepsilon=\frac{1}{2}\min\Bigl(\{1\}\cup \{\rho(A)-r:\mathbf{0}\in A\subseteq E,\ \rho(A)>r\}\Bigr)>0.$$
Every such $A$ then satisfies either $\rho(A)\le r$ or $\rho(A)\ge r+2\varepsilon$.

For the centers in the statement, put $A=\{\boldsymbol\lambda_i-\boldsymbol\lambda_0:0\le i\leq k\}$. Since $\varepsilon\le1/2$, we have
\begin{align*}
  |\boldsymbol\lambda_i-\boldsymbol\lambda_0|_1
  &\leq|\boldsymbol\lambda_i-\mathbf{x}_i|_1 +|\mathbf{x}_i-\mathbf{x}_0|_1 +|\mathbf{x}_0-\boldsymbol\lambda_0|_1\\
  &<2r+\varepsilon<2(r+1).
\end{align*}
Hence $\mathbf{0}\in A\subseteq E$. Using $\mathbf{x}_0-\boldsymbol\lambda_0$ as a candidate center gives
$$\rho(A)\le\max_i|\mathbf{x}_0-\boldsymbol\lambda_i|_1 \le\max_i\bigl(|\mathbf{x}_0-\mathbf{x}_i|_1 +|\mathbf{x}_i-\boldsymbol\lambda_i|_1\bigr) <r+\varepsilon.$$
By the choice of $\varepsilon$, this forces $\rho(A)\le r$. Thus $A$ is co-covering, and translating back proves the assertion.
\end{proof}

\subsection{Oriented chains and the covering construction}

For points $\mathbf{p}_0,\ldots,\mathbf{p}_k\in\R^3$, the symbol $[\mathbf{p}_0,\ldots,\mathbf{p}_k]$ denotes a {\it formal oriented $k$-simplex}. Permuting the vertices multiplies the symbol by the sign of the permutation, and a repeated vertex makes it zero. Affine dependence is allowed: $k$ specifies the number of vertices minus one, not the dimension of their convex hull. For $k=3$ we also call the symbol a {\it formal oriented tetrahedron}. Four distinct coplanar points therefore give a formal tetrahedron with zero three-dimensional volume but with its four oriented triangles still present in its boundary. Likewise, a list of five points in $\R^3$ defines a formal oriented $4$-simplex.

A finite integer linear combination of formal oriented $k$-simplices is an {\it integral $k$-chain}. Its boundary is defined linearly by
$$\partial[\mathbf{p}_0,\ldots,\mathbf{p}_k] =\sum_{i=0}^k(-1)^i [\mathbf{p}_0,\ldots,\widehat{\mathbf{p}_i},\ldots,\mathbf{p}_k] \qquad(k\geq 1),$$
where the hat means omission; a formal $0$-simplex has zero boundary. In a {\it periodic chain}, we identify
$$[\mathbf{p}_0+\mathbf{u},\ldots,\mathbf{p}_k+\mathbf{u}] =[\mathbf{p}_0,\ldots,\mathbf{p}_k] \qquad(\mathbf{u}\in L).$$
The same vector $\mathbf{u}$ translates the entire list; the vertices are not taken modulo $L$ separately. We combine coefficients with the signs prescribed by their vertex orders. Thus $\partial z=0$ means that every oriented face, up to common lattice translation, has total coefficient zero.

For example, take independent $\mathbf{a},\mathbf{b}\in L$ and put
$$t=[\mathbf{0},\mathbf{a},\mathbf{a}+\mathbf{b}] +[\mathbf{0},\mathbf{a}+\mathbf{b},\mathbf{b}].$$
The common diagonal cancels in the boundary, leaving
$$\partial t=[\mathbf{a},\mathbf{a}+\mathbf{b}]-[\mathbf{0},\mathbf{b}] +[\mathbf{a}+\mathbf{b},\mathbf{b}]+[\mathbf{0},\mathbf{a}].$$
In the periodic chain group,
$$[\mathbf{a},\mathbf{a}+\mathbf{b}]=[\mathbf{0},\mathbf{b}], \qquad [\mathbf{a}+\mathbf{b},\mathbf{b}]=[\mathbf{a},\mathbf{0}] =-[\mathbf{0},\mathbf{a}],$$
so $\partial t=0$. The geometric parallelogram has a boundary, but its opposite edges cancel after the periodic identification. The same cancellation for triangular faces gives the tetrahedron equations used in Section~\ref{sec:manual-topology}.

For a formal oriented tetrahedron $T=[\mathbf{p}_0,\mathbf{p}_1,\mathbf{p}_2,\mathbf{p}_3]$ with arbitrary real vertices, put
$$\mathcal D(T) =\det(\mathbf{p}_1-\mathbf{p}_0,\mathbf{p}_2-\mathbf{p}_0, \mathbf{p}_3-\mathbf{p}_0).$$
Its {\it signed volume} is $\mathcal D(T)/6$, and its ordinary volume is $|\mathcal D(T)|/6$. Extend $\mathcal D$ linearly to integral $3$-chains. Since $\mathcal D(T)$ is unchanged by common translation, it is also well defined on periodic chains.

Fix a basis matrix $M$ of $L$ with $\det M=D>0$. For a formal oriented tetrahedron $T=[\boldsymbol\ell_0,\boldsymbol\ell_1,\boldsymbol\ell_2,\boldsymbol\ell_3]$ with lattice vertices, define its {\it normalized determinant} by
\begin{equation*}
  \begin{split}
    d_T&=\frac{\mathcal D(T)}{D}\\
    &=\det\bigl(M^{-1}(\boldsymbol\ell_1-\boldsymbol\ell_0),
    M^{-1}(\boldsymbol\ell_2-\boldsymbol\ell_0),
    M^{-1}(\boldsymbol\ell_3-\boldsymbol\ell_0)\bigr)\in\Z.
  \end{split}
\end{equation*}
For a periodic chain $z=\sum_T z_TT$ of lattice tetrahedra, put
$$d(z)=\sum_T z_Td_T.$$
Here $T$ ranges over oriented common translation classes and $z_T\in\Z$ is its coefficient. The three quantities are related by
$$\mathcal D(z)=D\,d(z),\qquad \text{total signed volume of }z =\frac{\mathcal D(z)}{6}=\frac{D}{6}\,d(z).$$

\begin{proposition}\label{prop:nerve-degree}
Every covering $K+L=\R^3$ admits a finite integral chain $z$ of co-covering formal oriented tetrahedra, modulo common lattice translation, such that
\begin{equation}\label{eq:degree-six}
  \partial z=0,\qquad d(z)=6.
\end{equation}
\end{proposition}
\begin{proof}
We triangulate space into small tetrahedra and replace each vertex $\mathbf{x}$ by a lattice point $\boldsymbol\lambda(\mathbf{x})$ with $\mathbf{x}\in K+\boldsymbol\lambda(\mathbf{x})$. Lemma~\ref{lem:intersection-thickening} ensures that the resulting tetrahedra are co-covering. We will check that this replacement preserves both face cancellation and total signed volume.

Choose $\varepsilon$ as in Lemma~\ref{lem:intersection-thickening} and a positive integer $N$ such that
$$\frac{|M\mathbf{e}_1|_1+|M\mathbf{e}_2|_1+ |M\mathbf{e}_3|_1}{N}<\varepsilon.$$
Let
$$G=\frac{1}{N} L=\frac{1}{N} M\Z^3$$
be the grid, and call its points {\it grid vertices}. For each $\mathbf{a}\in\Z^3$, subdivide the small parallelepiped $M(\mathbf{a}+[0,1]^3)/N$ into the six tetrahedra with vertices
$$\frac{M\mathbf{a}}{N},\quad \frac{M(\mathbf{a}+\mathbf{e}_i)}{N},\quad \frac{M(\mathbf{a}+\mathbf{e}_i+\mathbf{e}_j)}{N},\quad \frac{M(\mathbf{a}+\mathbf{e}_1+\mathbf{e}_2+\mathbf{e}_3)}{N}, \qquad(i,j,k)\in\mathfrak{S}_3.$$
This is the Kuhn triangulation used in Section~\ref{sec:proof-outline}. If two small parallelepipeds share a parallelogram face with vertices
$$\mathbf{u},\quad \mathbf{u}+\frac{M\mathbf{e}_i}{N},\quad \mathbf{u}+\frac{M\mathbf{e}_j}{N},\quad \mathbf{u}+\frac{M(\mathbf{e}_i+\mathbf{e}_j)}{N},$$
both subdivisions use the diagonal from $\mathbf{u}$ to $\mathbf{u}+M(\mathbf{e}_i+\mathbf{e}_j)/N$. The subdivisions therefore agree on shared faces and give an $L$-periodic triangulation of space. No new vertices are introduced, so its vertex set is $G$. The choice of $N$ makes the diameter of each tetrahedron less than $\varepsilon$.

Order $G$ lexicographically by $M$-coordinates: the first coordinate at which two points differ determines their order. This order is preserved by lattice translation. Write the $6N^3$ tetrahedra in $M[0,1]^3$ as $T_\alpha$ with their vertices in this order, and choose $\sigma_\alpha\in\{1,-1\}$ so that $\sigma_\alpha T_\alpha$ has positive orientation. Put
$$c=\sum_{\alpha=1}^{6N^3}\sigma_\alpha T_\alpha.$$
Every triangular face not contained in the boundary of $M[0,1]^3$ belongs to two selected tetrahedra with opposite induced orientations, so these internal faces cancel. On the boundary of $M[0,1]^3$, matching triangles on opposite parallelogram faces differ by a lattice translation and again have opposite induced orientations. Consequently,
$$\partial c=0\qquad\text{in the periodic chain group}.$$
To choose the covering centers consistently under lattice translation, first take the $N^3$ representatives of $G/L$ given by
$$\mathbf{x}=\frac{M\mathbf{a}}{N}, \qquad \mathbf{a}\in\{0,\ldots,N-1\}^3.$$
For each representative, choose $\boldsymbol\lambda(\mathbf{x})\in L$ such that $\mathbf{x}\in K+\boldsymbol\lambda(\mathbf{x})$, and extend by
$$\boldsymbol\lambda(\mathbf{x}+\mathbf{u}) =\boldsymbol\lambda(\mathbf{x})+\mathbf{u} \qquad(\mathbf{u}\in L).$$
Denote vertex replacement by $\boldsymbol\lambda_*$. For $T_\alpha=[\mathbf{x}_{\alpha0},\ldots,\mathbf{x}_{\alpha3}]$, set
\begin{align*}
  T'_\alpha&=[\boldsymbol\lambda(\mathbf{x}_{\alpha0}),
  \boldsymbol\lambda(\mathbf{x}_{\alpha1}),
  \boldsymbol\lambda(\mathbf{x}_{\alpha2}),
  \boldsymbol\lambda(\mathbf{x}_{\alpha3})],\\
  z&=\sum_{\alpha=1}^{6N^3}\sigma_\alpha T'_\alpha
  =\boldsymbol\lambda_*c.
\end{align*}
A repeated vertex makes a term zero. In the alternating boundary formula, deleting a vertex before replacement gives the same face as deleting its image afterwards. The translation rule also preserves the identification of translated tetrahedra. Thus $\partial z=\boldsymbol\lambda_*(\partial c)=0$. Figure~\ref{fig:periodic-cancellation} illustrates these cancellations in two dimensions. Since each original tetrahedron has diameter less than $\varepsilon$, Lemma~\ref{lem:intersection-thickening} shows that every nonzero term of $z$ is co-covering.

\begin{figure}[tbp]
 \centering
\begin{tikzpicture}[x=1cm,y=1cm,line cap=round,line join=round,
                    font=\footnotesize,>=stealth]
  \node[font=\small] at (2,3.55) {(a) Periodic edges};
  \node[font=\small] at (9.4,3.55) {(b) Vertex replacement};

  \coordinate (o) at (.35,.6);
  \coordinate (a) at (3,.6);
  \coordinate (b) at (1,2.75);
  \coordinate (ab) at (3.65,2.75);
  \fill[black!4] (o)--(a)--(ab)--cycle;
  \fill[black!10] (o)--(ab)--(b)--cycle;
  \draw[line width=.55pt] (o)--(a)--(ab)--(b)--cycle;
  \draw[black!55,line width=.45pt] (o)--(ab);

  \draw[->,line width=.7pt] (1.15,.6)--(2.05,.6);
  \draw[->,line width=.7pt] (2.85,2.75)--(1.95,2.75);
  \draw[->,line width=.7pt] (3.2,1.26)--(3.43,2.02);
  \draw[->,line width=.7pt] (.82,2.15)--(.58,1.35);
  \draw[->,line width=.7pt] (2.70,2.01)--(1.38,1.15);
  \draw[->,line width=.7pt] (1.25,1.31)--(2.57,2.17);
  \foreach \point in {o,a,b,ab}
    \fill (\point) circle (1.25pt);
  \node[below left,inner sep=2pt] at (o) {$\mathbf{0}$};
  \node[below right,inner sep=2pt] at (a) {$\mathbf{a}$};
  \node[above left,inner sep=2pt] at (b) {$\mathbf{b}$};
  \node[above right,inner sep=2pt] at (ab) {$\mathbf{a}+\mathbf{b}$};
  \node at (2,-.05) {$\partial t=0$\quad periodically};

  \coordinate (x0) at (5.15,1.6);
  \coordinate (x1) at (7.15,1.6);
  \coordinate (x2) at (6.15,2.8);
  \coordinate (x3) at (6.15,.4);
  \filldraw[fill=black!4,draw=black!65,line width=.5pt]
    (x0)--(x1)--(x2)--cycle;
  \filldraw[fill=black!12,draw=black!65,line width=.5pt]
    (x1)--(x0)--(x3)--cycle;
  \draw[->,line width=.75pt] (5.65,1.67)--(6.65,1.67);
  \draw[->,line width=.75pt] (6.65,1.53)--(5.65,1.53);
  \foreach \point in {x0,x1,x2,x3}
    \fill (\point) circle (1.25pt);
  \node[left,inner sep=2pt] at (x0) {$\mathbf{x}_0$};
  \node[right,inner sep=2pt] at (x1) {$\mathbf{x}_1$};
  \node[above,inner sep=2pt] at (x2) {$\mathbf{x}_2$};
  \node[below,inner sep=2pt] at (x3) {$\mathbf{x}_3$};
  \node at (6.15,-.25)
    {$[\mathbf{x}_0,\mathbf{x}_1]+[\mathbf{x}_1,\mathbf{x}_0]=0$};

  \draw[->,line width=.65pt] (7.9,1.6)--(9.2,1.6)
    node[midway,above,inner sep=4pt,font=\normalsize] {$\boldsymbol\lambda_*$};

  \coordinate (y0) at (10,.6);
  \coordinate (y1) at (12.7,.6);
  \coordinate (y2) at (11.65,2.9);
  \coordinate (y3) at (11.35,1.9);
  \filldraw[fill=black!4,draw=black!65,line width=.5pt]
    (y0)--(y1)--(y2)--cycle;
  \filldraw[fill=black!12,draw=black!65,line width=.5pt]
    (y1)--(y0)--(y3)--cycle;
  \draw[->,line width=.75pt] (10.8,.67)--(11.9,.67);
  \draw[->,line width=.75pt] (11.9,.53)--(10.8,.53);
  \foreach \point in {y0,y1,y2,y3}
    \fill (\point) circle (1.25pt);
  \node[below left,inner sep=2pt] at (y0) {$\mathbf{y}_0$};
  \node[below right,inner sep=2pt] at (y1) {$\mathbf{y}_1$};
  \node[above,inner sep=2pt] at (y2) {$\mathbf{y}_2$};
  \node[above left,inner sep=2pt] at (y3) {$\mathbf{y}_3$};
  \node at (11.35,-.25)
    {$[\mathbf{y}_0,\mathbf{y}_1]+[\mathbf{y}_1,\mathbf{y}_0]=0$};
\end{tikzpicture}
  \caption{Boundary cancellation in two dimensions. Left: the diagonal cancels internally, and opposite boundary edges cancel after common lattice translation. Right: with $\mathbf{y}_i=\boldsymbol\lambda(\mathbf{x}_i)$, the shared edge has opposite orientations in the two summands before and after vertex replacement, even when their images overlap.}
 \label{fig:periodic-cancellation}
\end{figure}

It remains to show that the replacement preserves total signed volume. We do this by constructing a formal $4$-chain $H$ with $z-c=\partial H$ and proving that $\mathcal D(\partial H)=0$. The prism identity~\cite[proof of Theorem~2.10]{Hatcher2002} gives $H=Pc$, where $P$ is defined as follows. For a simplex $T=[\mathbf{x}_0,\ldots,\mathbf{x}_k]$ of the grid triangulation in its fixed lexicographic order, put
$$PT=\sum_{i=0}^k(-1)^i [\mathbf{x}_0,\ldots,\mathbf{x}_i, \boldsymbol\lambda(\mathbf{x}_i),\ldots,\boldsymbol\lambda(\mathbf{x}_k)].$$
For other orientations use the corresponding sign, and extend linearly. The summands are formal oriented simplices with real vertices. Every face inherits the same lexicographic order. Since lattice translation preserves this order and the rule for $\boldsymbol\lambda$, the definition of $P$ also respects periodic identification.

To check the identity, put $\mathbf{y}_i=\boldsymbol\lambda(\mathbf{x}_i)$ and
$$R_i=[\mathbf{x}_0,\ldots,\mathbf{x}_{i-1}, \mathbf{y}_i,\ldots,\mathbf{y}_k] \quad(0\le i\le k+1),$$
with $R_0=[\mathbf{y}_0,\ldots,\mathbf{y}_k]$ and $R_{k+1}=[\mathbf{x}_0,\ldots,\mathbf{x}_k]$. In the $i$th signed summand of $PT$, deleting $\mathbf{x}_i$ gives $R_i$, while deleting $\mathbf{y}_i$ gives $-R_{i+1}$. The remaining faces sum to $-P(\partial T)$. Therefore
\begin{align*}
  \partial(PT)
  &=\sum_{i=0}^k(R_i-R_{i+1})-P(\partial T)\\
  &=\boldsymbol\lambda_*T-T-P(\partial T).
\end{align*}
Equivalently,
$$\partial P+P\partial=\boldsymbol\lambda_*-\mathrm{id}.$$
Since $\partial c=0$, its application to $c$ gives
\begin{equation}\label{eq:prism-chain-boundary}
  z-c=\partial(Pc).
\end{equation}

For a formal oriented $4$-simplex $S=[\mathbf{p}_0,\ldots,\mathbf{p}_4]$, write $\mathbf{p}_j=(p_{1j},p_{2j},p_{3j})^{\mathsf T}$. Let $B_j$ be the $4\times4$ matrix with columns $(1,\mathbf{p}_\ell^{\mathsf T})^{\mathsf T}$ for $0\leq \ell\leq 4$, $\ell\neq j$, in increasing order. Expansion along the first row gives
\begin{align*}
  0
  &=\det\begin{pmatrix}
    1&1&1&1&1\\
    1&1&1&1&1\\
    p_{10}&p_{11}&p_{12}&p_{13}&p_{14}\\
    p_{20}&p_{21}&p_{22}&p_{23}&p_{24}\\
    p_{30}&p_{31}&p_{32}&p_{33}&p_{34}
  \end{pmatrix}\\
  &=\sum_{j=0}^4(-1)^j\det B_j\\
  &=\sum_{j=0}^4(-1)^j
  \mathcal D[\mathbf{p}_0,\ldots,\widehat{\mathbf{p}_j},\ldots,\mathbf{p}_4]
  =\mathcal D(\partial S).
\end{align*}
The determinant is zero because its first two rows are equal. Subtracting the first column of each $B_j$ from its other columns gives the corresponding tetrahedron determinant in the third line. By linearity, $\mathcal D$ vanishes on the boundary of every formal $4$-chain. Applying it to~\eqref{eq:prism-chain-boundary} yields
$$\mathcal D(z)-\mathcal D(c)=\mathcal D(\partial Pc)=0.$$
Each positively oriented grid tetrahedron has determinant $D/N^3$, so
$$\mathcal D(c)=6N^3\frac{D}{N^3}=6D, \qquad d(z)=\frac{\mathcal D(z)}{D}=6.$$
\end{proof}

The chain may have negative coefficients and overlapping tetrahedra. The identity $d(z)=6$ specifies its total signed volume and does not assert a tiling. Proposition~\ref{prop:nerve-degree} provides two properties needed below. Since $d(z)\neq 0$, the chain contains a nondegenerate tetrahedron. The identity $\partial z=0$ gives relations among the coefficients of tetrahedra sharing a triangular face up to lattice translation.

From now on, assume $D>4$. Every edge vector of a co-covering tetrahedron belongs to $K-K=3O$. The next two sections establish the determinant and parity bounds for these vectors that will be used in Section~\ref{sec:active-core}.

\section{A determinant bound for four vectors}\label{sec:four-vector}

We first prove a determinant estimate for three vectors and their sum. For a lattice tetrahedron with vertices $\mathbf p_0,\ldots,\mathbf p_3$, the successive edge vectors satisfy
$$(\mathbf p_1-\mathbf p_0)+(\mathbf p_2-\mathbf p_1)+(\mathbf p_3-\mathbf p_2)=\mathbf p_3-\mathbf p_0.$$
If every edge vector has $\ell_1$ norm at most three, the three summands and their sum satisfy the same norm bound. Under $D>4$, the estimate below will imply that the tetrahedron is unimodular.

\begin{lemma}\label{lem:four-vector}
Let $B$ be a real $3\times3$ matrix, and let $\mathbf{a},\mathbf{b},\mathbf{c}\in\R^3$ be linearly independent. If
$$|B\mathbf{a}|_1,\ |B\mathbf{b}|_1,\ |B\mathbf{c}|_1, \ |B(\mathbf{a}+\mathbf{b}+\mathbf{c})|_1\leq 3,$$
then
\begin{equation}\label{eq:four-vector-bound}
  |\det B|\leq\frac{8}{|\det[\mathbf{a}\ \mathbf{b}\ \mathbf{c}]|}.
\end{equation}
The constant $8$ is best possible.
\end{lemma}

\begin{proof}
The four vectors and their negatives are the images of the eight vertices of $[-1,1]^3$ under the matrix
$$Q=\tfrac{1}{2} B[\mathbf{b}+\mathbf{c}\ \mathbf{a}+\mathbf{c}\ \mathbf{a}+\mathbf{b}].$$
Then $Q(1,1,1)^{\mathsf T}=B(\mathbf{a}+\mathbf{b}+\mathbf{c})$ and $Q(-1,1,1)^{\mathsf T}=B\mathbf{a}$. The other sign choices show that the eight images of the cube vertices are precisely $\pm B\mathbf{a},\pm B\mathbf{b},\pm B\mathbf{c}$ and $\pm B(\mathbf{a}+\mathbf{b}+\mathbf{c})$. In particular, all eight have $\ell_1$ norm at most three. The transformation from the columns $B\mathbf{a},B\mathbf{b},B\mathbf{c}$ to those of $Q$ has determinant $1/4$, so
\begin{equation}\label{eq:four-vector-change}
  |\det Q|=\tfrac{1}{4}|\det B|\,|\det[\mathbf{a}\ \mathbf{b}\ \mathbf{c}]|.
\end{equation}
For a row vector $\mathbf{x}\in\R^3$, take the average of $|\mathbf{x}\cdot\boldsymbol\varepsilon|$ over $\boldsymbol\varepsilon\in\{\pm1\}^3$:
$$\nu(\mathbf{x}):=\tfrac{1}{8}\sum_{\boldsymbol\varepsilon\in\{\pm1\}^3} |\mathbf{x}\cdot\boldsymbol\varepsilon|.$$
If $\mathbf{q}_1,\mathbf{q}_2,\mathbf{q}_3$ are the rows of $Q$, then
$$\sum_{i=1}^3\nu(\mathbf{q}_i) =\tfrac{1}{8}\sum_{\boldsymbol\varepsilon\in\{\pm1\}^3} |Q\boldsymbol\varepsilon|_1\leq 3.$$
We shall bound $|\det Q|$ in terms of these three row averages.

The average is unchanged by signs and permutations of the coordinates. Writing their absolute values as $\alpha\ge\beta\ge\gamma\ge0$ and pairing opposite cube vertices gives
\begin{align*}
  \nu(\mathbf{x})
  &=\tfrac{1}{4}\bigl((\alpha+\beta+\gamma)+(\alpha+\beta-\gamma)
  +(\alpha-\beta+\gamma)+|\alpha-\beta-\gamma|\bigr)\\
  &=\max\{\alpha,\tfrac{1}{2}(\alpha+\beta+\gamma)\}
  =\max\{|\mathbf{x}|_\infty,\tfrac{1}{2}|\mathbf{x}|_1\}.
\end{align*}
Thus $\nu$ is a norm with unit ball
$$P=\{\mathbf{x}:\nu(\mathbf{x})\leq 1\} =[-1,1]^3\cap\{\mathbf{x}:|\mathbf{x}|_1\leq 2\}.$$
Each cube face $x_i=\pm1$ is cut down to the square $|x_j|+|x_k|\le1$, where $\{i,j,k\}=\{1,2,3\}$; the remaining faces are triangles. Hence the vertices of $P$ are the twelve permutations of $(\pm1,\pm1,0)$.

We first show that a matrix with all three rows in $P$ has determinant of absolute value at most two. Fix two rows $\mathbf{y},\mathbf{z}$ and write the third as a convex combination $\mathbf{x}=\sum_jt_j\mathbf{v}_j$ of vertices of $P$. Linearity and the triangle inequality give
$$|\det[\mathbf{x};\mathbf{y};\mathbf{z}]| \leq \sum_jt_j|\det[\mathbf{v}_j;\mathbf{y};\mathbf{z}]| \le\max_j|\det[\mathbf{v}_j;\mathbf{y};\mathbf{z}]|.$$
Applying this successively to the three rows reduces the determinant bound to matrices $[\mathbf{v}_1;\mathbf{v}_2;\mathbf{v}_3]$ whose rows are vertices of $P$. Their determinants are integers, and each row has Euclidean norm $\sqrt2$. Hadamard's inequality gives
$$|\det[\mathbf{v}_1;\mathbf{v}_2;\mathbf{v}_3]| \leq (\sqrt2)^3=2\sqrt2<3,$$
so every such determinant has absolute value at most two. For arbitrary nonzero rows $\mathbf{x}_i$, apply this bound to $\mathbf{x}_i/\nu(\mathbf{x}_i)$ and multiply by $\prod_i\nu(\mathbf{x}_i)$. If a row is zero, the determinant is zero. In either case,
$$|\det[\mathbf{x}_1;\mathbf{x}_2;\mathbf{x}_3]| \leq 2\prod_{i=1}^3 \nu(\mathbf{x}_i).$$
Applying this to $Q$ and using the arithmetic--geometric mean inequality, we obtain
$$|\det Q|\leq 2\prod_{i=1}^3\nu(\mathbf{q}_i) \leq 2\left(\frac{\sum_{i=1}^3\nu(\mathbf{q}_i)}{3}\right)^3\le2.$$
This proves the claim by \eqref{eq:four-vector-change}. Equality is attained at
$$
  Q=\begin{pmatrix}1&1&0\\1&-1&0\\0&0&1\end{pmatrix}:
  \qquad |\det Q|=2,\qquad |Q\boldsymbol\varepsilon|_1=3.
$$
The invertible change defining $Q$ gives equality in \eqref{eq:four-vector-bound}.
\end{proof}

\begin{corollary}\label{cor:four-vector-cut}
Suppose $L=B\Z^3$ and $|\det B|>4$. There do not exist $\mathbf{a},\mathbf{b},\mathbf{c}\in\Z^3$ such that
$$|B\mathbf{a}|_1,\ |B\mathbf{b}|_1,\ |B\mathbf{c}|_1, \ |B(\mathbf{a}+\mathbf{b}+\mathbf{c})|_1\leq 3, \qquad |\det[\mathbf{a}\ \mathbf{b}\ \mathbf{c}]|\geq 2.$$
\end{corollary}
A nondegenerate lattice tetrahedron is {\it unimodular} if the three edge vectors at a vertex form a basis of the lattice. This condition is independent of the chosen vertex.

\begin{corollary}\label{cor:short-unimodular}
Let $L=B\Z^3$ with $|\det B|>4$. Every nondegenerate lattice tetrahedron whose edge vectors have $\ell_1$ norm at most three is unimodular.
\end{corollary}

\begin{proof}
After a lattice translation, we may assume that the nondegenerate lattice tetrahedron has vertices $\mathbf{0}, B\mathbf{a}, B\mathbf{b}, B\mathbf{c},$ where $\mathbf{a},\mathbf{b},\mathbf{c}\in\Z^3$. The three successive edge vectors have integer coordinates $\mathbf{a},\mathbf{b}-\mathbf{a},\mathbf{c}-\mathbf{b}$, and their sum is $\mathbf{c}$. The images of these four vectors under $B$ are edge vectors of the tetrahedron and therefore have $\ell_1$ norm at most three. Lemma~\ref{lem:four-vector} gives
$$|\det B|\,|\det(\mathbf{a},\mathbf{b},\mathbf{c})| =|\det B|\,|\det(\mathbf{a},\mathbf{b}-\mathbf{a},\mathbf{c}-\mathbf{b})|\leq 8.$$
The second factor is a nonzero integer, and $|\det B|>4$, so it equals one. Hence the three edge vectors at zero form a basis of $L$.
\end{proof}

Corollary~\ref{cor:short-unimodular} supplies a lattice basis from every nondegenerate tetrahedron in the covering chain. Relative to such a basis, the determinant restriction in Corollary~\ref{cor:four-vector-cut} will also bound the coordinates of a fourth point joined to the vertices of a triangular face. Lemma~\ref{lem:degenerate-bridge} uses these bounds in both the coplanar and noncoplanar cases.

\section{Short lattice vectors and parity}\label{sec:short-vectors}

Throughout this section, assume
$$K+L=\R^3,\qquad K=\tfrac32O,\qquad D=\det L>4.$$
We prove two restrictions on short lattice vectors. The lower bound $|\mathbf v|_1>3/2$ for nonzero lattice vectors implies that every short vector is primitive and has nonzero parity. We then show that each nonzero class of $L/2L$ contains at most two short directions. These restrictions will be used to compare edge directions and to exclude an additional apex in Lemma~\ref{lem:face-stability}.

\subsection{An overlap estimate}

Call $\mathbf{v}\in L\setminus\{\mathbf{0}\}$ {\it short} if $|\mathbf{v}|_1\le3$, and identify a {\it direction} with the pair $\{\mathbf{v},-\mathbf{v}\}$. A lattice vector is {\it primitive} if it is not an integer multiple $k\mathbf{w}$, with $\mathbf{w}\in L$ and $|k|\ge2$.

\begin{lemma}\label{lem:short-minimum}
For every $\mathbf{v}\in L\setminus\{\mathbf{0}\}$,
\begin{equation}\label{eq:one-overlap}
  D+\operatorname{vol}\bigl(K\cap(K+\mathbf{v})\bigr)\le\frac{9}{2}.
\end{equation}
Consequently,
\begin{equation}\label{eq:short-minimum}
  \lambda_1(L):=\min_{\mathbf{v}\in L\setminus\{\mathbf{0}\}}|\mathbf{v}|_1
  >\frac{3}{2}.
\end{equation}
In particular, every short vector is primitive, and no nonzero short vector belongs to $2L$.
\end{lemma}

Here $\lambda_1(L)$ is the first successive minimum of $L$ with respect to the unit octahedron $O$.

\begin{proof}
Fix $\mathbf{v}\in L\setminus\{\mathbf{0}\}$. For $\mathbf{x}\in\R^3$, put
$$C_{\mathbf{x}}=\{\boldsymbol\ell\in L:\mathbf{x}\in K+\boldsymbol\ell\}, \qquad n(\mathbf{x})=|C_{\mathbf{x}}|,$$
and let
$$e_{\mathbf{v}}(\mathbf{x}) =\#\{\boldsymbol\ell\in C_{\mathbf{x}}: \boldsymbol\ell+\mathbf{v}\in C_{\mathbf{x}}\}.$$
The set $C_{\mathbf{x}}$ is finite and nonempty. It cannot be closed under addition of $\mathbf{v}$, since repeated addition would then produce infinitely many points. Thus
$$e_{\mathbf{v}}(\mathbf{x})\leq n(\mathbf{x})-1.$$
Let $\Omega$ be a half-open fundamental parallelepiped of $L$. Since the sets $\Omega-\boldsymbol\ell$, $\boldsymbol\ell\in L$, partition $\R^3$, every nonnegative integrable function $f$ satisfies
\begin{align*}
  \int_\Omega\sum_{\boldsymbol\ell\in L}
  f(\mathbf{x}-\boldsymbol\ell)\,d\mathbf{x}
  &=\sum_{\boldsymbol\ell\in L}
  \int_{\Omega-\boldsymbol\ell}f(\mathbf{y})\,d\mathbf{y}\\
  &=\int_{\R^3}f(\mathbf{y})\,d\mathbf{y}.
\end{align*}
For $f=\mathbf{1}_K$, the sum in the first integral is $n(\mathbf{x})$. For $f=\mathbf{1}_{K\cap(K+\mathbf{v})}$, it is $e_{\mathbf{v}}(\mathbf{x})$. Consequently,
\begin{align*}
  D+\operatorname{vol}\bigl(K\cap(K+\mathbf{v})\bigr)
  &=\int_\Omega(1+e_{\mathbf{v}}(\mathbf{x}))\,d\mathbf{x}\\
  &\leq \int_\Omega n(\mathbf{x})\,d\mathbf{x}
  =\operatorname{vol}K=\frac{9}{2},
\end{align*}
which proves \eqref{eq:one-overlap}.

To prove the strict lower bound, suppose, for a contradiction, that $|\mathbf{v}|_1\leq 3/2$. For every $\mathbf{x}\in\mathbf{v}/2+(3/4)O$, the triangle inequality gives
$$\max\{|\mathbf{x}|_1,|\mathbf{x}-\mathbf{v}|_1\} \leq \frac{|\mathbf{v}|_1}{2} +\left|\mathbf{x}-\frac{\mathbf{v}}{2}\right|_1 \leq\frac{3}{2}.$$
Thus
$$\frac{\mathbf{v}}{2}+\frac{3}{4}O\subset K\cap(K+\mathbf{v}).$$
The set on the left has volume $(4/3)(3/4)^3=9/16$. Combining this with~\eqref{eq:one-overlap}, we obtain
$$D\le\frac{9}{2}-\frac{9}{16}=\frac{63}{16}<4,$$
contrary to $D>4$. Thus~\eqref{eq:short-minimum} holds. A nonprimitive short vector $\mathbf{v}=k\mathbf{w}$ would give $|\mathbf{w}|_1\le3/2$, which is impossible. The same argument excludes a nonzero short vector in $2L$.
\end{proof}

Fix a basis matrix $B$ of $L$ with $\det B=D>0$. For $\mathbf{v}\in\Z^3\setminus\{\mathbf{0}\}$, let $[\mathbf{v}]$ be the member of $\{\mathbf{v},-\mathbf{v}\}$ whose first nonzero coordinate is positive, and define
\begin{equation}\label{eq:short-direction-set}
  S=\{[\mathbf v]:\mathbf v\in\Z^3\setminus\{\mathbf0\},|B\mathbf v|_1\leq3\},\qquad \pm S=S\cup(-S).
\end{equation}
Thus $\pm S$ consists of all nonzero integer vectors $\mathbf v$ satisfying $|B\mathbf v|_1\leq 3$. These vectors are primitive by Lemma~\ref{lem:short-minimum}.

\subsection{At most two directions in a parity class}\label{sec:parity-matrix}

Two lattice vectors have the same {\it parity} if their difference belongs to $2L$. The quotient $L/2L$ has eight classes, one of which is $2L$ itself. Opposite vectors have the same parity, so parity is also defined for a short direction. For an integer coordinate vector $\mathbf v\in\Z^3$, its parity is its reduction $\overline{\mathbf{v}}\in\mathbb F_2^3$, where $\mathbb F_2=\Z/2\Z$. Lemma~\ref{lem:short-minimum} implies that $\overline{\mathbf v}\neq \mathbf0$ whenever $\mathbf v\in S$. The following lemma bounds the number of short directions in each nonzero parity class.

\begin{lemma}\label{lem:parity-two}
Each nonzero class in $L/2L$ contains at most two short directions.
\end{lemma}

To prove Lemma~\ref{lem:parity-two}, we will assume that three distinct short directions have the same parity and construct a matrix $A$ for which $D=2|\det A|$. Its columns have bounded norms, and the covering condition provides an additional bound on a linear combination of these columns. The following estimate uses these constraints to force $D<4$, contradicting our standing assumption, so we establish it first.

\begin{lemma}\label{lem:parity-matrix}
Let $C=(c_{ij})$ be a real $3\times3$ matrix whose columns have $\ell_1$ norm at most one. If there is $\mathbf{m}\in\R^3$ such that
$$|m_j|\geq 1\quad(1\leq j\leq 3),\qquad |C\mathbf{m}|_1\leq 2,$$
then
$$|\det C|\leq \frac{125}{216}.$$
\end{lemma}

\begin{proof}
The conclusion holds if $|\det C|\le1/2$, so assume $|\det C|>1/2$. Write $\mathbf{c}_j$ for column $j$ of $C$. Expand the determinant into its six signed products. Each product contains one entry from each column, in three distinct rows. Bounding the entry from column one by $|\mathbf{c}_1|_\infty$ gives
\begin{align*}
  |\det C|&\leq|\mathbf{c}_1|_\infty\sum_{j\neq k}|c_{j2}|\,|c_{k3}|\\
  &\leq|\mathbf{c}_1|_\infty \left(\sum_j|c_{j2}|\right)\left(\sum_k|c_{k3}|\right)
  \leq|\mathbf{c}_1|_\infty.
\end{align*}
The same argument with either of the other columns gives $|\det C|\leq|\mathbf{c}_j|_\infty$ for every $j$. Thus every column has an entry of absolute value greater than $1/2$; it is unique because the column has $\ell_1$ norm at most one.

The three entries just selected occupy different rows. Indeed, suppose that the large entries in columns one and two occupy the same row, with absolute values $a,b>1/2$. Every two-by-two minor of these columns using that row would have absolute value at most
$$a(1-b)+b(1-a)=\frac{1}{2}-2(a-1/2)(b-1/2)<\frac{1}{2}.$$
A minor not using that row is bounded by $(1-a)(1-b)<1/4$. Let $\Delta_i$ be the two-by-two determinant obtained by deleting row $i$ from the first two columns. The preceding bounds give $|\Delta_i|<1/2$ for all $i$. Expansion along column three therefore gives
$$|\det C|=\left|\sum_{i=1}^3(-1)^{i+3}c_{i3}\Delta_i\right| \leq\max_i|\Delta_i|\sum_{i=1}^3|c_{i3}| <\tfrac{1}{2},$$
contrary to our assumption.

For each $j$, replace $\mathbf{c}_j$ by $\operatorname{sgn}(m_j)\mathbf{c}_j$ and $m_j$ by $|m_j|$. This leaves $C\mathbf{m}$ unchanged and gives $m_j\ge1$. Permute the rows and change their signs to put the three large entries on a positive diagonal. These operations preserve the hypotheses and $|\det C|$. We have strict column dominance:
\begin{equation}\label{eq:strict-column-dominance}
  c_{jj}>\sum_{i\ne j}|c_{ij}|\qquad(1\leq j\leq 3).
\end{equation}
For $0\leq t\leq 1$, let $C(t)$ have the same diagonal as $C$ and off-diagonal entries $t c_{ij}$. Each $C(t)$ is invertible. Indeed, if $C(t)^{\mathsf T}\mathbf{x}=\mathbf{0}$ with $\mathbf{x}\ne\mathbf{0}$, choose $k$ with $|x_k|=\max_i|x_i|>0$. The $k$th equation gives
$$c_{kk}|x_k|\le t\sum_{i\neq k}|c_{ik}|\,|x_i| \leq \left(\sum_{i\neq k}|c_{ik}|\right)|x_k|<c_{kk}|x_k|,$$
a contradiction. Hence $\det C(t)$ never vanishes. Since $\det C(0)=c_{11}c_{22}c_{33}>0$ and the determinant is continuous, $\det C=\det C(1)>0$.

Put $s_j=\sum_i|c_{ij}|\le1$ and $p_j=\sum_i c_{ij}>0$. Positivity of $p_j$ follows from~\eqref{eq:strict-column-dominance}. The original hypothesis gives
\begin{equation}\label{eq:parity-column-sums}
  \sum_jp_j\leq \sum_jm_jp_j
  =\sum_i(C\mathbf{m})_i\leq |C\mathbf{m}|_1\leq 2.
\end{equation}

We now use only strict column dominance, $s_j\le1$, and the bound $\sum_jp_j\le2$ in~\eqref{eq:parity-column-sums}. To increase the determinant while preserving $s_j$ and $p_j$, we move each positive off-diagonal entry to the diagonal of its column. If $q=c_{ij}>0$ with $i\ne j$, replace $c_{ij}$ by zero and $c_{jj}$ by $c_{jj}+q$. This preserves both $s_j$ and $p_j$: in column $j$ it subtracts $q$ from one positive entry and adds $q$ to another. It also increases $c_{jj}-\sum_{i\ne j}|c_{ij}|$ by $2q$, so \eqref{eq:strict-column-dominance} continues to hold. Write $\operatorname{cof}_{ij}(C)$ for $(-1)^{i+j}$ times the determinant of the matrix obtained by deleting row $i$ and column $j$ from $C$. Expansion along column $j$ shows that the change in determinant is
$$q\bigl(\operatorname{cof}_{jj}(C)-\operatorname{cof}_{ij}(C)\bigr)>0.$$
For the strict inequality, simultaneous relabeling of rows and columns reduces to $j=1,i=2$, where
\begin{align*}
  \operatorname{cof}_{11}(C)-\operatorname{cof}_{21}(C)
  &=c_{33}(c_{22}+c_{12})-c_{32}(c_{13}+c_{23})\\
  &\geq c_{33}(c_{22}-|c_{12}|)-|c_{32}|(|c_{13}|+|c_{23}|)\\
  &>|c_{32}|\bigl(c_{33}-|c_{13}|-|c_{23}|\bigr)\geq 0.
\end{align*}
Here the strict inequality follows from $c_{22}>|c_{12}|+|c_{32}|$, and the last inequality from $c_{33}>|c_{13}|+|c_{23}|$. The cofactors in column $j$ depend on the other two columns, which the replacement leaves unchanged. Thus the displayed determinant change is exact. The determinant-sign and cofactor comparisons used here are standard properties of diagonally dominant matrices~\cite{DK2011}.

After all positive off-diagonal entries have been removed, the matrix $H$ has positive diagonal and nonpositive off-diagonal entries. Since the column norms and sums were preserved,
$$h_{jj}=\frac{s_j+p_j}{2},\qquad \sum_jh_{jj}\le\frac{3+2}{2}=\frac{5}{2}.$$
In the determinant expansion, the three transposition terms and the two three-cycle terms are nonpositive. Thus only the diagonal product can be positive, and
$$|\det C|=\det C\leq \det H \le h_{11}h_{22}h_{33} \leq \left(\frac{h_{11}+h_{22}+h_{33}}{3}\right)^3 \leq \left(\frac{5}{6}\right)^3.$$
\end{proof}

\begin{proof}[Proof of Lemma~\ref{lem:parity-two}]
Suppose that a nonzero parity class contains three distinct short directions, represented by lattice vectors $\mathbf{u},\mathbf{v},\mathbf{w}\in L$, and put
$$\mathbf{a}=\mathbf{u}/2,\qquad \mathbf{b}=\mathbf{v}/2,\qquad \mathbf{c}=\mathbf{w}/2.$$
We first show that $\mathbf{a},\mathbf{b},\mathbf{c}$ are linearly independent. The points $\pm\mathbf{a},\pm\mathbf{b},\pm\mathbf{c}$ are distinct because $\mathbf{u},\mathbf{v},\mathbf{w}$ represent distinct nonzero directions. Their norms are at most $3/2$, and their pairwise differences belong to $L$: for example, $\mathbf{a}-\mathbf{b}=(\mathbf{u}-\mathbf{v})/2\in L$ and $\mathbf{a}-(-\mathbf{a})=\mathbf{u}\in L$. Thus all six lie in $K$ and in one coset of $L$. More generally, points of one lattice coset contained in $K$ are distinct modulo $2L$: a difference $2\boldsymbol\ell\neq \mathbf{0}$ would imply $|\boldsymbol\ell|_1\leq 3/2$, contradicting Lemma~\ref{lem:short-minimum}.

Put $W=\operatorname{span}\{\mathbf{a},\mathbf{b},\mathbf{c}\}$ and $\Gamma=L\cap W$. Translating the six points by $-\mathbf{a}$ puts them in $\Gamma$ without changing their pairwise differences. Since $2\Gamma\subseteq2L$, they are also distinct modulo $2\Gamma$. Hence
$$6\leq |\Gamma/2\Gamma|=2^{\operatorname{rank}\Gamma} \leq 2^{\dim W}\le8.$$
It follows that $\dim W=3$, so $\mathbf{a},\mathbf{b},\mathbf{c}$ are linearly independent.

Let $A=[\mathbf{a}\ \mathbf{b}\ \mathbf{c}]$ and
$$D_3=\{\mathbf{z}\in\mathbb Z^3:z_1+z_2+z_3\text{ is even}\}.$$
The lattice vectors $-2\mathbf{a},\mathbf{b}-\mathbf{a},\mathbf{c}-\mathbf{a}$ generate $AD_3$: their $A$-coordinates are $-2\mathbf{e}_1,\mathbf{e}_2-\mathbf{e}_1,\mathbf{e}_3-\mathbf{e}_1$, which generate the index-two subgroup $D_3$ of $\Z^3$. Thus $AD_3\subset L$. Put $k=[L:AD_3]$, a positive integer since $AD_3$ is a full-rank sublattice of $L$. A sublattice of index $k$ has $k$ times the determinant, while $A$ multiplies volume by $|\det A|$. Since $\det D_3=2$, the column norm bound gives
$$kD=\det(AD_3)=2|\det A| \leq 2|\mathbf{a}|_1|\mathbf{b}|_1|\mathbf{c}|_1\leq \frac{27}{4}.$$
Since $D>4$, we have $k<27/16<2$, so $k=1$. Hence
$$L=AD_3,\qquad D=2|\det A|.$$

It remains to use Lemma~\ref{lem:parity-matrix} to bound this determinant. Cover $A(1/2,1/2,1/2)^{\mathsf T}$ by a translate centered at $A\mathbf{z}$, where $\mathbf{z}\in D_3$. The vector $\mathbf{m}=(1,1,1)^{\mathsf T}-2\mathbf{z}$ has odd integer coordinates and satisfies
$$|A\mathbf{m}|_1\leq 3.$$
The columns of $\widetilde A=(2/3)A$ have $\ell_1$ norm at most one, and $|\widetilde A\mathbf{m}|_1\le2$. Applying Lemma~\ref{lem:parity-matrix} yields
$$D=2\left(\frac{3}{2}\right)^3|\det\widetilde A| \leq \frac{125}{32}<4,$$
contrary to the standing assumption. Thus three short directions cannot occur.
\end{proof}

\section{A common set of eight edge directions}\label{sec:active-core}

Starting from Proposition~\ref{prop:nerve-degree}, we will obtain a chain of nondegenerate co-covering tetrahedra with cancelling faces and nonzero determinant sum. We then choose lattice coordinates in which all its edge vectors belong to the signed set represented in~\eqref{eq:square-core}. This will reduce the face equations in Section~\ref{sec:manual-topology} to a fixed finite list of tetrahedron classes.

We continue under $D>4$ and use integer coordinates relative to the basis matrix $B$ fixed in Section~\ref{sec:short-vectors}. A set of integer points is co-covering when its image under $B$ lies in a translate of $K$. A tetrahedron with integer vertices is called a {\it short-edge tetrahedron} if every difference between two distinct vertices belongs to $\pm S$.

\subsection{A closed chain of nondegenerate tetrahedra}

Since determinants of integer vectors are integers, Corollary~\ref{cor:four-vector-cut} gives
\begin{equation}\label{eq:core-four-vectors}
  \mathbf x,\mathbf y,\mathbf z,\mathbf x+\mathbf y+\mathbf z \in\pm S
  \quad\Longrightarrow\quad
  |\det(\mathbf x,\mathbf y,\mathbf z)|\leq 1.
\end{equation}
We use this bound to restrict the coordinates of a point joined by short vectors to the vertices of a triangular face. The coplanar case will identify the degenerate terms that can share a face with a nondegenerate tetrahedron. The noncoplanar case will be used in Lemma~\ref{lem:shared-face-apices}.

\begin{lemma}\label{lem:degenerate-bridge}
Let $\{\mathbf0,\mathbf p,\mathbf q,\mathbf r\}$ be a nondegenerate short-edge tetrahedron. Suppose that $\mathbf x\in\Z^3$ satisfies
$$\mathbf x,\quad \mathbf x-\mathbf p,\quad \mathbf x-\mathbf q\in\pm S.$$
Write
$$\mathbf x=m\mathbf p+n\mathbf q+h\mathbf r.$$
Then $m,n,h\in\Z$ and
$$|m|\leq 1,\qquad |n|\leq 1,\qquad |h|\leq 1,\qquad |m+n+h-1|\leq 1.$$
In particular, every co-covering, affinely dependent set of four distinct integer points sharing a triangular face, up to common lattice translation, with a nondegenerate short-edge tetrahedron is the vertex set of a parallelogram.
\end{lemma}

\begin{proof}
By Corollary~\ref{cor:short-unimodular}, $(\mathbf p,\mathbf q,\mathbf r)$ is a basis of $\Z^3$, so $m,n,h$ are integers. Each summand and each right-hand side in the following identities belongs to $\pm S$.
\begin{align*}
  \mathbf x+(\mathbf p-\mathbf x)+(\mathbf r-\mathbf p)
  &=\mathbf r,\\
  \mathbf x+(\mathbf q-\mathbf x)+(\mathbf r-\mathbf q)
  &=\mathbf r,\\
  (\mathbf x-\mathbf p)+(\mathbf q-\mathbf x)
  +(\mathbf r-\mathbf q)
  &=\mathbf r-\mathbf p,\\
  \mathbf p+(\mathbf q-\mathbf p)+(\mathbf x-\mathbf q)
  &=\mathbf x.
\end{align*}
Since $|\det(\mathbf p,\mathbf q,\mathbf r)|=1$, the triples of summands on the left have absolute determinants
$$|n|,\qquad |m|,\qquad |m+n+h-1|,\qquad |h|,$$
respectively. Equation~\eqref{eq:core-four-vectors} gives the required bounds.

For the final assertion, translate the shared face to $\{\mathbf0,\mathbf p,\mathbf q\}$ and denote the apex of the nondegenerate tetrahedron by $\mathbf r$. Write the coplanar set as
$$P=\{\mathbf0,\mathbf p,\mathbf q,\mathbf x\}.$$
Since $P$ is co-covering and $K-K=3O$, its three differences $\mathbf x,\mathbf x-\mathbf p, \mathbf x-\mathbf q$ belong to $\pm S$. Affine dependence gives $h=0$, so the bounds above imply
$$m,n\in\{-1,0,1\},\qquad 0\leq m+n\leq 2.$$
The four points are distinct, so $\mathbf x\notin\{\mathbf0,\mathbf p,\mathbf q\}$. Consequently,
$$\mathbf x\in \{\mathbf p+\mathbf q,\ \mathbf p-\mathbf q,\mathbf q-\mathbf p\}.$$
In each case the four points can be partitioned into two pairs with the same sum. Since $\mathbf p,\mathbf q$ are linearly independent, $P$ is the vertex set of a nondegenerate parallelogram.
\end{proof}

The parallelogram conclusion will allow us to separate the nondegenerate tetrahedra while preserving face cancellation. We will also retain a face-connected part with nonzero determinant sum. We first specify the adjacency used in this statement. Following the periodic identification in Section~\ref{sec:lattice-topology}, two triangular faces $F,F'$ belong to the same {\it unoriented face class} if
$$F'=F+\boldsymbol\ell \qquad\text{for some }\boldsymbol\ell\in\Z^3,$$
where the faces are regarded as vertex sets. For a finite family of nondegenerate tetrahedra, form a graph whose vertices are their unoriented four-point translation classes. Join two distinct vertices when the corresponding tetrahedra share an unoriented triangular face class. We call the family {\it face-connected} when this graph is connected. Orientations affect the chain coefficients but are ignored in this graph. Combine the coefficients of tetrahedra in each translation class according to their orientations. The {\it support} consists of the unoriented classes with nonzero total coefficients.

\begin{proposition}\label{prop:closed-component}
There is a finite integral periodic chain $c$ of nondegenerate co-covering tetrahedra whose support is face-connected and satisfies $\partial c=0$ and $d(c)\neq 0$.
\end{proposition}

\begin{proof}
Let $z$ be the chain from Proposition~\ref{prop:nerve-degree}, expressed in integer coordinates relative to $B$. Since $\det B=D>0$, the determinant in these coordinates equals the normalized determinant $d$. Thus
$$\partial z=0,\qquad d(z)=6.$$

Let $\iota$ denote reflection through the origin, acting on an oriented tetrahedron by
$$\iota[\mathbf v_0,\mathbf v_1,\mathbf v_2,\mathbf v_3]= [-\mathbf v_0,-\mathbf v_1,-\mathbf v_2,-\mathbf v_3].$$
Apply this reflection to every tetrahedron of $z$, keeping its coefficient, and put $w=z-\iota z$. Reflection commutes with taking boundaries and changes the sign of each three-dimensional determinant. Hence
$$\partial w=\partial z-\iota(\partial z)=0, \qquad d(w)=d(z)-d(\iota z)=2d(z)=12.$$
Since $K=-K$, every reflected tetrahedron is still co-covering.

Let
$$T=[\mathbf0,\mathbf a,\mathbf b,\mathbf a+\mathbf b]$$
be an oriented parallelogram. Reflecting its vertices and then translating by $\mathbf a+\mathbf b$ gives $[\mathbf a+\mathbf b,\mathbf b,\mathbf a,\mathbf0]$. Two exchanges of vertices restore the original order, so the orientation sign is unchanged. Hence, in the periodic chain group,
$$\iota T=[\mathbf a+\mathbf b,\mathbf b,\mathbf a,\mathbf0]=T.$$
Every parallelogram class therefore has coefficient zero in $w=z-\iota z$.

After combining coefficients, let $w_{\rm nd}$ consist of the nondegenerate tetrahedra in $w$. By the final assertion of Lemma~\ref{lem:degenerate-bridge}, a degenerate co-covering tetrahedron sharing a triangular face with a nondegenerate short-edge tetrahedron must be a parallelogram. The same conclusion holds when the faces agree only after lattice translation. Since all parallelogram classes have been removed, no remaining degenerate tetrahedron shares a face class with $w_{\rm nd}$. The cancellation of these faces in $\partial w=0$ therefore gives $\partial w_{\rm nd}=0$. Degenerate tetrahedra have determinant zero, so
$$\partial w_{\rm nd}=0,\qquad d(w_{\rm nd})=12.$$
Partition the support of $w_{\rm nd}$ into the connected components of its face graph, and let $c_1,\ldots,c_m$ be the corresponding chains. Different components share no triangular face class. Thus face cancellation takes place within each component, and
$$\partial c_i=0\quad(1\leq i\leq m), \qquad \sum_{i=1}^{m}d(c_i)=d(w_{\rm nd})=12.$$
At least one $c_i$ has $d(c_i)\neq 0$. Taking $c=c_i$ gives the required finite integral periodic chain of nondegenerate co-covering tetrahedra with face-connected support.
\end{proof}

\subsection{Comparing tetrahedra across a face}

The chain in Proposition~\ref{prop:closed-component} has face-connected support. We now apply the coordinate bounds of Lemma~\ref{lem:degenerate-bridge} when the fourth point lies outside the shared face plane. This gives the possible apices and a fourth common edge direction. We then classify the union of two different tetrahedron direction sets having at least four directions in common.

\begin{lemma}\label{lem:shared-face-apices}
Suppose two nondegenerate short-edge tetrahedra share the face $\{\mathbf0,\mathbf p,\mathbf q\}$, and let $\mathbf r$ be the apex of the first tetrahedron. Then the apex of the second tetrahedron belongs to
\begin{equation*}
  \begin{split}
    \{\mathbf r,\ \mathbf r\pm\mathbf p,\ \mathbf r\pm\mathbf q, \ \mathbf r\pm(\mathbf p-\mathbf q), \ -\mathbf r+\mathbf p,\ -\mathbf r+\mathbf q, \ -\mathbf r+\mathbf p+\mathbf q\}.
  \end{split}
\end{equation*}
In particular, the two tetrahedra have at least four common edge directions.
\end{lemma}

\begin{proof}
By Corollary~\ref{cor:short-unimodular}, $\mathbf p,\mathbf q,\mathbf r$ form a basis of $\Z^3$. We may therefore write the second apex as
$$\mathbf x=m\mathbf p+n\mathbf q+h\mathbf r.$$
The edges joining $\mathbf x$ to the shared face give $\mathbf x,\mathbf x-\mathbf p,\mathbf x-\mathbf q\in\pm S$. Lemma~\ref{lem:degenerate-bridge} now yields
$$m,n,h\in\Z,\qquad |m|,|n|,|h|\leq 1, \qquad |m+n+h-1|\leq 1.$$
The second tetrahedron is nondegenerate, so $h\neq 0$. Hence $h=1$ or $h=-1$.

If $h=1$, the last bound becomes $|m+n|\leq 1$. Together with $m,n\in\{-1,0,1\}$, this gives
$$\mathbf x-\mathbf r\in \{\mathbf0,\pm\mathbf p,\pm\mathbf q,\pm(\mathbf p-\mathbf q)\}.$$
If $h=-1$, the bounds give $1\le m+n\leq 2$, so
$$(m,n)\in\{(1,0),(0,1),(1,1)\}.$$
These are exactly the stated ten candidates.

We next show that the two tetrahedra share an edge direction outside the plane of their common face. For $h=1$, we may write
$$\mathbf x=\mathbf r+\mathbf u-\mathbf v,\qquad \mathbf u,\mathbf v\in\{\mathbf0,\mathbf p,\mathbf q\}.$$
Consequently,
$$\mathbf x-\mathbf u=\mathbf r-\mathbf v.$$
For $h=-1$, we may write
$$\mathbf x=-\mathbf r+\mathbf u+\mathbf v, \qquad \mathbf u,\mathbf v\in\{\mathbf0,\mathbf p,\mathbf q\}, \quad \mathbf u\ne\mathbf v.$$
Then
$$\mathbf x-\mathbf u=-(\mathbf r-\mathbf v).$$
In each case these are edge vectors of the two tetrahedra with the same direction. Since $\mathbf r-\mathbf v\notin \operatorname{span}\{\mathbf p,\mathbf q\}$, they give a fourth common edge direction.
\end{proof}

To compare tetrahedra with four common edge directions, we will use the following implication.
\begin{equation}\label{eq:parallel-segment-obstruction}
  \mathbf u\pm\mathbf v,\ \mathbf w\pm\mathbf v\in\pm S,
  \quad \mathbf u,\mathbf v,\mathbf w\in\Z^3
  \quad\Longrightarrow\quad
  \det(\mathbf u,\mathbf v,\mathbf w)=0.
\end{equation}
In the identity
$$(\mathbf w+\mathbf v)-(\mathbf u+\mathbf v)+(\mathbf u-\mathbf v)=\mathbf w-\mathbf v,$$
the three summands and their sum all belong to $\pm S$. Therefore, applying~\eqref{eq:core-four-vectors} gives
$$2|\det(\mathbf u,\mathbf v,\mathbf w)|=|\det(\mathbf w+\mathbf v,-\mathbf u-\mathbf v, \mathbf u-\mathbf v)|\leq 1.$$
Since $\det(\mathbf u,\mathbf v,\mathbf w)$ is an integer, it must be zero.

For a lattice basis $\mathbf{p},\mathbf{q},\mathbf{r}$ of $\Z^3$, write
\begin{equation}\label{eq:square-core}
  C=\{\mathbf{p},\ \mathbf{q},\ \mathbf{p}+\mathbf{q},\ \mathbf{p}-\mathbf{q}, \ \mathbf{r},\ \mathbf{r}+\mathbf{p}, \ \mathbf{r}+\mathbf{q}, \ \mathbf{r}+\mathbf{p}+\mathbf{q}\}.
\end{equation}
The first four vectors lie in $\operatorname{span}\{\mathbf p,\mathbf q\}$. The remaining four form the vertex set of a fundamental parallelogram in $\mathbf r+\operatorname{span}\{\mathbf p,\mathbf q\}$. The next lemma shows that this set represents every eight-direction union of two tetrahedra with at least four common edge directions. The tetrahedra need not share a face.

\begin{lemma}\label{lem:two-tetrahedra}
Suppose two nondegenerate short-edge tetrahedra have different edge-direction sets with at least four directions in common. Then the union of their direction sets has seven or eight directions. In the seven-direction case, all seven nonzero parity classes occur exactly once. In the eight-direction case, the union is represented by~\eqref{eq:square-core} for a suitable basis of $\Z^3$.
\end{lemma}

\begin{proof}
Each tetrahedron has six edge directions, so the hypotheses already imply that the union has seven or eight directions. We determine its form in each case.

Choose four common directions. The corresponding four edges of a tetrahedron form either a triangle with one further edge or a four-cycle. Three edges of a nondegenerate tetrahedron have linearly dependent vectors exactly when they form a triangle. Thus, if the four edges form a four-cycle in one tetrahedron, they also form a four-cycle in the other, because any three of their vectors are independent.

In the four-cycle case, direct the edges consecutively around the first cycle and write their vectors as $\mathbf u,\mathbf v,\mathbf w,\mathbf t$. Then
$$\mathbf u+\mathbf v+\mathbf w+\mathbf t=\mathbf0.$$
Since any three are independent, every linear relation among the four vectors is a scalar multiple of this one. Corresponding edge vectors in the second cycle agree up to sign, because they are primitive. Their sum is also zero, so these signs must either all remain unchanged or all be reversed. The two remaining edges of each tetrahedron are the diagonals of its four-cycle. Each diagonal vector is the sum of two consecutive cycle vectors. Since complementary pairs have opposite sums, the two diagonal directions in either tetrahedron form a two-element subset of
$$[\mathbf u+\mathbf v],\qquad [\mathbf u+\mathbf w],\qquad [\mathbf u+\mathbf t].$$
The two pairs therefore share a direction, giving a fifth common edge direction. Five edges of a tetrahedron contain a triangle and a further edge. Consequently, in both cases the tetrahedra have three common directions forming a triangular face in each, together with another common edge direction.

A triangle with primitive integer edge vectors is determined by its three edge directions up to translation and central reflection. We may therefore translate and, if necessary, centrally reflect the second tetrahedron so that the two triangular faces coincide. These operations preserve edge directions and the condition that edge vectors belong to $\pm S$.

Choose a fourth common edge direction, distinct from the three edge directions of the shared face. In the first tetrahedron, the corresponding edge joins the apex to a face vertex. Translate both tetrahedra so that this face vertex is the origin. Writing the other two face vertices as $\mathbf p,\mathbf q$ and the first apex as $\mathbf r$, the first tetrahedron has vertices
$$\{\mathbf0,\mathbf p,\mathbf q,\mathbf r\},$$
and $[\mathbf r]$ is the chosen fourth common direction. By Corollary~\ref{cor:short-unimodular}, $\mathbf p,\mathbf q,\mathbf r$ form a basis of $\Z^3$. In the second tetrahedron, this direction also joins its apex to one of $\mathbf0,\mathbf p,\mathbf q$, with displacement $\mathbf r$ or $-\mathbf r$. Hence the second apex belongs to
$$\{\mathbf r,-\mathbf r, \mathbf p+\mathbf r,\mathbf p-\mathbf r, \mathbf q+\mathbf r,\mathbf q-\mathbf r\}.$$

The choice $\mathbf r$ gives the same direction set and is excluded by hypothesis. The choice $-\mathbf r$ gives
$$\mathbf p\pm\mathbf r,\qquad \mathbf q\pm\mathbf r\in\pm S.$$
Equation~\eqref{eq:parallel-segment-obstruction} would then imply $\det(\mathbf p,\mathbf r,\mathbf q)=0$, contrary to the independence of the basis vectors.

If the second apex is $\mathbf p-\mathbf r$, its three edge vectors to the shared face are
$$\mathbf p-\mathbf r,\qquad -\mathbf r,\qquad \mathbf p-\mathbf q-\mathbf r.$$
The first two directions already occur in the first tetrahedron, and the last is new. Thus the union has seven directions. The six original directions are represented by
$$\mathbf p,\quad\mathbf q,\quad\mathbf r,\quad \mathbf p-\mathbf q,\quad \mathbf p-\mathbf r,\quad \mathbf q-\mathbf r.$$
Since $\mathbf p,\mathbf q,\mathbf r$ form a lattice basis, these occupy six distinct nonzero parity classes, missing only $\overline{\mathbf p}+\overline{\mathbf q}+\overline{\mathbf r}$. The new vector satisfies
$$\overline{\mathbf p-\mathbf q-\mathbf r}=\overline{\mathbf p}+\overline{\mathbf q} +\overline{\mathbf r},$$
so all seven nonzero classes occur exactly once. The apex $\mathbf q-\mathbf r$ gives the same conclusion after interchanging $\mathbf p$ and $\mathbf q$.

The remaining apices are $\mathbf p+\mathbf r$ and $\mathbf q+\mathbf r$, which are also exchanged by $\mathbf p\leftrightarrow\mathbf q$. For the apex $\mathbf p+\mathbf r$, the three edge vectors to the shared face are
$$\mathbf p+\mathbf r,\qquad \mathbf r,\qquad \mathbf p+\mathbf r-\mathbf q.$$
The first and third give two new directions. The union therefore has eight directions, represented by
$$\mathbf p,\quad\mathbf q,\quad\mathbf r,\quad \mathbf p-\mathbf q,\quad\mathbf p-\mathbf r,\quad \mathbf q-\mathbf r,\quad \mathbf p+\mathbf r,\quad\mathbf p+\mathbf r-\mathbf q.$$
The vectors in~\eqref{eq:square-core} for the basis $(-\mathbf p,-\mathbf r,\mathbf q)$ represent exactly these eight directions.
\end{proof}

\begin{figure}[htbp]
\centering
\begin{tikzpicture}[x=1cm,y=1cm,line cap=round,line join=round,
                    font=\small]
  \coordinate (o)  at (0,0);
  \coordinate (p)  at (2.7,-.3);
  \coordinate (q)  at (1,1);
  \coordinate (pq) at (3.7,.7);
  \coordinate (pmq) at (1.7,-1.3);
  \coordinate (r)  at (-.5,2.15);
  \coordinate (rp) at (2.2,1.85);
  \coordinate (rq) at (.5,3.15);
  \coordinate (rpq) at (3.2,2.85);

  \fill[black!3] (o)--(p)--(pq)--(q)--cycle;
  \draw[black!40,line width=.35pt] (p)--(pq)--(q);
  \draw[black!35,densely dashed,line width=.35pt]
    (o)--(r) (p)--(rp) (q)--(rq) (pq)--(rpq);
  \filldraw[fill=black!5,draw=black!70,line width=.6pt]
    (r)--(rp)--(rpq)--(rq)--cycle;

  \draw[->,>=stealth,line width=.7pt] (o)--(p);
  \draw[->,>=stealth,line width=.7pt] (o)--(q);
  \draw[->,>=stealth,line width=.7pt] (o)--(pq);
  \draw[->,>=stealth,line width=.7pt] (o)--(pmq);
  \fill (o) circle (1.25pt);
  \foreach \point in {r,rp,rq,rpq}
    \fill (\point) circle (1.6pt);

  \node[anchor=east,inner sep=3pt] at (o) {$\mathbf{0}$};
  \node[anchor=north west,inner sep=3pt] at (p) {$\mathbf{p}$};
  \node[anchor=south east,inner sep=3pt] at (q) {$\mathbf{q}$};
  \node[anchor=west,inner sep=3pt] at (pq) {$\mathbf{p}+\mathbf{q}$};
  \node[anchor=north,inner sep=3pt] at (pmq) {$\mathbf{p}-\mathbf{q}$};
  \node[anchor=east,inner sep=3pt] at (r) {$\mathbf{r}$};
  \node[anchor=north west,inner sep=3pt] at (rp) {$\mathbf{r}+\mathbf{p}$};
  \node[anchor=south east,inner sep=3pt] at (rq) {$\mathbf{r}+\mathbf{q}$};
  \node[anchor=south west,inner sep=3pt] at (rpq)
    {$\mathbf{r}+\mathbf{p}+\mathbf{q}$};
\end{tikzpicture}
  \caption{An affine schematic of~\eqref{eq:square-core} in basis coordinates. The lower arrows represent four planar directions; the four upper points form a fundamental parallelogram. Dashed segments indicate translation by $\mathbf{r}$.}
\label{fig:eight-directions}
\end{figure}

\subsection{Extension across a shared face}

We next show that $\pm C$ continues to contain all edge vectors when we pass to another nondegenerate short-edge tetrahedron sharing a face. This requires all eight directions of $C$ to be short, as in the eight-direction case of Lemma~\ref{lem:two-tetrahedra}.

\begin{lemma}\label{lem:face-stability}
Suppose all eight directions in~\eqref{eq:square-core} are short. If a nondegenerate short-edge tetrahedron shares a triangular face with a nondegenerate tetrahedron whose edges belong to $\pm C$, then all its edges also belong to $\pm C$.
\end{lemma}

\begin{proof}
The three edge vectors of the shared face already belong to $\pm C$. The other three edge vectors of the short-edge tetrahedron belong to $\pm S$, and it remains to show that they also belong to $\pm C$. For convenience, call the third coordinate in the basis $(\mathbf p,\mathbf q,\mathbf r)$ of $\Z^3$ the height, and put
$$H=\{\mathbf{p},\mathbf{q},\mathbf{p}+\mathbf{q},\mathbf{p}-\mathbf{q}\},\qquad Q=\mathbf{r}+\{\mathbf{0},\mathbf{p},\mathbf{q},\mathbf{p}+\mathbf{q}\}.$$
Every edge vector of the shared triangular face belongs to $\pm C$ and therefore has height $0$ or $\pm 1$. Translate both tetrahedra so that a vertex of minimum height on the shared face becomes $\mathbf 0$. Each vertex of this face then has height $0$ or $1$. If all three face vertices have the same height, the maximum and minimum values of each of their first two coordinates differ by at most one. Since these coordinates are integers, the face vertices correspond to three vertices of a translate of $[0,1]^2$. If the face vertices do not all have the same height, reflect and translate both tetrahedra if necessary so that two vertices of the shared face have height zero and the third has height one. Write the two height-zero vertices as $\mathbf0,\mathbf h$ and the height-one vertex as $\mathbf c$. Then $\mathbf h\in\pm H$ and $\mathbf c,\mathbf c-\mathbf h\in Q$. The edge joining $\mathbf0$ to $\mathbf h$ has vector $\mathbf h=\mathbf c-(\mathbf c-\mathbf h)$. The first two coordinates of $\mathbf h$ therefore give a side or diagonal vector of $[0,1]^2$.

We now choose one representative shared face for each of the three cases above. The following changes give another basis of $\Z^3$
\begin{align*}
  (\mathbf p,\mathbf q,\mathbf r)
  &\longmapsto(-\mathbf p,\mathbf q,\mathbf r+\mathbf p),\\
  (\mathbf p,\mathbf q,\mathbf r)
  &\longmapsto(\mathbf p,-\mathbf q,\mathbf r+\mathbf q),\\
  (\mathbf p,\mathbf q,\mathbf r)
  &\longmapsto(\mathbf q,\mathbf p,\mathbf r).
\end{align*}
Substituting any of these bases into the definitions of $H$ and $Q$ leaves the sets $\pm H$ and $Q$ unchanged. Thus the signed vector set $\pm C$ remains unchanged. We keep the notation $\mathbf p,\mathbf q,\mathbf r$ after changing the basis. We may also apply a common translation or a central reflection to both tetrahedra. Their face remains shared, and their edge vectors are unchanged up to sign. If all three face vertices have the same height, translate the vertex of the coordinate square adjacent to the other two face vertices to the origin. The basis changes then allow us to write the other two vertices as $\mathbf p,\mathbf q$. If the edge joining the two height-zero vertices corresponds to a side of $[0,1]^2$, the basis changes allow us to take $\mathbf h=\mathbf p$. The conditions $\mathbf c,\mathbf c-\mathbf p\in Q$ give
$$\mathbf c\in \{\mathbf r+\mathbf p,\mathbf r+\mathbf p+\mathbf q\}.$$
The second basis change turns the latter choice into $\mathbf r+\mathbf p$ in the new notation. If that edge corresponds to a diagonal of $[0,1]^2$, we may take $\mathbf h=\mathbf p+\mathbf q$. Then $\mathbf c,\mathbf c-\mathbf p-\mathbf q\in Q$ forces $\mathbf c=\mathbf r+\mathbf p+\mathbf q$. Thus the shared face can be written as one of the three vertex sets
$$F_0=\{\mathbf{0},\mathbf{p},\mathbf{q}\},\qquad F_1=\{\mathbf{0},\mathbf{p},\mathbf{r}+\mathbf{p}\},\qquad F_2=\{\mathbf{0},\mathbf{p}+\mathbf{q},\mathbf{r}+\mathbf{p}+\mathbf{q}\}.$$

Fix a shared face $F\in\{F_0,F_1,F_2\}$ and write $F=\{\mathbf0,\mathbf b,\mathbf c\}$. Let $\mathbf x$ be the fourth vertex of the short-edge tetrahedron whose edges we must show belong to $\pm C$. We first construct nondegenerate tetrahedra sharing $F$ whose edges already belong to $\pm C$. Their possible fourth vertices form the set
$$A_C(F)=\{\mathbf a\in\pm C:\mathbf a-\mathbf b,\mathbf a-\mathbf c\in\pm C,\ \det(\mathbf a,\mathbf b,\mathbf c)\neq 0\}.$$
For each $\mathbf a\in A_C(F)$, the tetrahedron $F\cup\{\mathbf a\}$ is nondegenerate and has all its edges in $\pm C\subseteq\pm S$. We use it as an auxiliary tetrahedron. Apply Lemma~\ref{lem:shared-face-apices} to $F\cup\{\mathbf a\}$ and $F\cup\{\mathbf x\}$, with $\mathbf a$ as the fixed apex. It gives
$$\mathbf x\in V(\mathbf a;\mathbf b,\mathbf c),$$
where the ten possible positions of $\mathbf x$ form the set
\begin{align*}
  V(\mathbf a;\mathbf b,\mathbf c)
  =\{&\mathbf a,\mathbf a\pm\mathbf b,\mathbf a\pm\mathbf c, \mathbf a\pm(\mathbf b-\mathbf c),-\mathbf a+\mathbf b,-\mathbf a+\mathbf c, -\mathbf a+\mathbf b+\mathbf c\}.
\end{align*}
The same vertex $\mathbf x$ must satisfy this condition for every $\mathbf a\in A_C(F)$, so it belongs to the intersection of these ten-point sets. Table~\ref{tab:face-apices} summarizes the possible positions of $\mathbf x$ allowed by these comparisons. The second column follows by applying the definition of $A_C(F)$ to $F=F_0,F_1,F_2$. If $\mathbf x\in A_C(F)$, all its edges to the shared face already belong to $\pm C$. It therefore remains to exclude the additional candidates
$$\bigcap_{\mathbf a\in A_C(F)}V(\mathbf a;\mathbf b,\mathbf c)\setminus A_C(F),$$
which are listed in the third column.

For example, the two auxiliary apices of $F_0$ give
$$
  \begin{aligned}
    V(-\mathbf{r};\mathbf{p},\mathbf{q})
    &\cap V(\mathbf{r}+\mathbf{p}+\mathbf{q};\mathbf{p},\mathbf{q})=\{-\mathbf{r},\ \mathbf{r}+\mathbf{p}+\mathbf{q},\
    -\mathbf{r}-\mathbf{p},\ -\mathbf{r}-\mathbf{q},\
    \mathbf{r}+\mathbf{p},\ \mathbf{r}+\mathbf{q}\}.
  \end{aligned}
$$
The first two are the apices in $A_C(F_0)$; the other four are exactly the extra candidates in its row of Table~\ref{tab:face-apices}.

\begin{table}[htbp]
\centering\small
  \caption{Auxiliary apices and the extra candidates left by their ten-point constraints.}
\label{tab:face-apices}
\begin{tabular}{@{}lll@{}}
\toprule
Face&$A_C(F)$&Extra candidates\\\midrule
$F_0$&$-\mathbf{r},\ \mathbf{r}+\mathbf{p}+\mathbf{q}$&$\begin{gathered}-\mathbf{r}-\mathbf{p},\ -\mathbf{r}-\mathbf{q},\ \mathbf{r}+\mathbf{p},\ \mathbf{r}+\mathbf{q}\end{gathered}$\\[4pt]
$F_1$&$-\mathbf{q},\ \mathbf{p}-\mathbf{q},\ \mathbf{r}+\mathbf{p}+\mathbf{q}$&$-\mathbf{r}-\mathbf{q}$\\
$F_2$&$\mathbf{p},\ \mathbf{q}$&$\begin{gathered}\mathbf{r}+\mathbf{p},\ \mathbf{r}+\mathbf{q},\ \mathbf{r}+2\mathbf{p}+\mathbf{q},\ \mathbf{r}+\mathbf{p}+2\mathbf{q}\end{gathered}$\\[4pt]
\bottomrule
\end{tabular}
\end{table}

We now exclude all positions of $\mathbf x$ listed in the third column of Table~\ref{tab:face-apices}. For $F_1$, suppose that $\mathbf x=-\mathbf r-\mathbf q$. The edge joining $\mathbf x$ to $\mathbf r+\mathbf p$ gives
$$\mathbf x-(\mathbf r+\mathbf p)=-(\mathbf p+\mathbf q+2\mathbf r)\in\pm S.$$
Together with $\mathbf p+\mathbf q,\mathbf p-\mathbf q\in C$, this yields three distinct short directions of the same nonzero parity, contrary to Lemma~\ref{lem:parity-two}.

For the remaining candidate apices for $F_0$ and $F_2$, choose $\mathbf u,\mathbf v,\mathbf w$ from the corresponding row of Table~\ref{tab:face-parallel-pairs}. With these choices, $\mathbf w+\mathbf v,\mathbf w-\mathbf v\in\pm C$, and each of $\mathbf u+\mathbf v,\mathbf u-\mathbf v$ belongs to $\pm C$ or is an edge vector of the short-edge tetrahedron $F\cup\{\mathbf x\}$. Thus all four vectors belong to $\pm S$, and \eqref{eq:parallel-segment-obstruction} gives
$$\det(\mathbf u,\mathbf v,\mathbf w)=0.$$
However, in every row, $\mathbf u$ has height $1$ or $-1$, and $\mathbf v,\mathbf w$ are $\mathbf p,\mathbf q$ in some order. Hence
$$|\det(\mathbf u,\mathbf v,\mathbf w)|=|\det(\mathbf p,\mathbf q,\mathbf r)|=1,$$
a contradiction.

\begin{table}[htbp]
\centering\small
  \caption{Choices of $\mathbf u,\mathbf v,\mathbf w$ used in~\eqref{eq:parallel-segment-obstruction} to rule out the remaining positions of $\mathbf x$.}
\label{tab:face-parallel-pairs}
\begin{tabular}{@{}ccccc@{}}
\toprule
Face&Apex&$\mathbf{u}$&$\mathbf{v}$&$\mathbf{w}$\\\midrule
$F_0$&$-\mathbf{r}-\mathbf{p}$&$-\mathbf{r}-\mathbf{p}$&$\mathbf{p}$&$\mathbf{q}$\\
$F_0$&$-\mathbf{r}-\mathbf{q}$&$-\mathbf{r}-\mathbf{q}$&$\mathbf{q}$&$\mathbf{p}$\\
$F_0$&$\mathbf{r}+\mathbf{p}$&$\mathbf{r}+\mathbf{p}$&$\mathbf{q}$&$\mathbf{p}$\\
$F_0$&$\mathbf{r}+\mathbf{q}$&$\mathbf{r}+\mathbf{q}$&$\mathbf{p}$&$\mathbf{q}$\\
$F_2$&$\mathbf{r}+\mathbf{p}$&$\mathbf{r}$&$\mathbf{q}$&$\mathbf{p}$\\
$F_2$&$\mathbf{r}+\mathbf{q}$&$\mathbf{r}$&$\mathbf{p}$&$\mathbf{q}$\\
$F_2$&$\mathbf{r}+2\mathbf{p}+\mathbf{q}$&$\mathbf{r}+\mathbf{p}+\mathbf{q}$&$\mathbf{p}$&$\mathbf{q}$\\
$F_2$&$\mathbf{r}+\mathbf{p}+2\mathbf{q}$&$\mathbf{r}+\mathbf{p}+\mathbf{q}$&$\mathbf{q}$&$\mathbf{p}$\\
\bottomrule
\end{tabular}
\end{table}

Hence $\mathbf x\in A_C(F)$, so all six edge vectors lie in $\pm C$.
\end{proof}

\subsection{Face-connected families of tetrahedra}

We show that all edge vectors of a face-connected family of nondegenerate short-edge tetrahedra belong to the same set $\pm C$ for a suitable lattice basis. If at least eight edge directions occur, we first find two tetrahedra whose union has exactly eight directions. Lemma~\ref{lem:two-tetrahedra} gives $C$ as in~\eqref{eq:square-core}, with all eight directions short. Applying Lemma~\ref{lem:face-stability} successively to tetrahedra sharing a triangular face then shows that all edge vectors of the family belong to this same $\pm C$. If only six or seven directions occur, we choose the basis directly to obtain the same containment. In this case, the unused directions of $C$ need not be short.

\begin{lemma}\label{lem:connected-eight-directions}
Every finite face-connected family of nondegenerate short-edge tetrahedra has all its edge vectors in $\pm C$ for one lattice basis $\mathbf{p},\mathbf{q},\mathbf{r}$ of $\Z^3$ and $C$ as in~\eqref{eq:square-core}.
\end{lemma}

\begin{proof}
Since the given family of tetrahedra is face-connected, we can label its members $T_1,\ldots,T_N$ so that each $T_i$ with $i>1$ shares a triangular face, up to lattice translation, with $T_j$ for some $j<i$. Lemma~\ref{lem:shared-face-apices} shows that any pair of these tetrahedra sharing a face has at most eight directions. If some such pair has an eight-direction union, Lemma~\ref{lem:two-tetrahedra} gives a basis of $\Z^3$ in which these directions are represented by $C$ in~\eqref{eq:square-core}. All eight directions of $C$ are short. Using face-connectedness again, we apply Lemma~\ref{lem:face-stability} successively across shared faces to show that all edge vectors of the family belong to this same $\pm C$.

We may therefore assume that every pair sharing a face has at most seven directions in its union. When $T_i$ is added, its earlier neighbor already has six directions, so $T_i$ introduces at most one direction not present among $T_1,\ldots,T_{i-1}$. We begin with $T_1$, which has six distinct edge directions. If a seventh direction first appears when $T_k$ is added, choose an earlier tetrahedron $T_j$ sharing a face with $T_k$. Before this addition, all directions are precisely the six directions of $T_j$. Thus the set $A$ of all directions present after adding $T_k$ is the seven-direction union of $T_j$ and $T_k$. By Lemma~\ref{lem:two-tetrahedra}, each nonzero parity class occurs exactly once in $A$.

Now, suppose that the whole family has more than seven edge directions. Let $T_\ell$ be the first tetrahedron to introduce an eighth direction $[\mathbf y]$. Before $T_\ell$ is added, the directions already present are exactly those in $A$. Since each nonzero parity class occurs exactly once in $A$, there is a unique $[\mathbf x]\in A$ with $\overline{\mathbf x}=\overline{\mathbf y}$. Choose an earlier tetrahedron $T_i$, with $i<\ell$, containing $[\mathbf x]$. Each of $T_i$ and $T_\ell$ has six edge directions, all belonging to the eight-element set $A\cup\{[\mathbf y]\}$. Thus, they have at least $4$ common directions. Their union contains the distinct same-parity directions $[\mathbf x]$ and $[\mathbf y]$, so the seven-direction case of Lemma~\ref{lem:two-tetrahedra} is impossible. That lemma gives a basis of $\Z^3$ for which their eight directions are represented by $C$ in~\eqref{eq:square-core}. The conclusion now follows from Lemma~\ref{lem:face-stability} as in the first paragraph.

It remains to consider families with six or seven edge directions. Suppose first that the family has seven edge directions. Their set is $A$, which was shown above to be the union of the edge-direction sets of two tetrahedra. The seven-direction case in the proof of Lemma~\ref{lem:two-tetrahedra} shows that, after translation and relabeling, one tetrahedron has vertices
$$\{\mathbf0,\mathbf a,\mathbf b,\mathbf c\},$$
and the additional direction is $[\mathbf a-\mathbf b-\mathbf c]$. By Corollary~\ref{cor:short-unimodular}, $\mathbf a,\mathbf b,\mathbf c$ form a basis of $\Z^3$. Set
$$\mathbf p=-\mathbf b,\qquad \mathbf q=-\mathbf c,\qquad \mathbf r=\mathbf a.$$
These also form a basis of $\Z^3$. The seven directions are then represented by
$$\mathbf p,\quad \mathbf q,\quad \mathbf r,\quad \mathbf p-\mathbf q,\quad \mathbf r+\mathbf p,\quad \mathbf r+\mathbf q,\quad \mathbf r+\mathbf p+\mathbf q.$$
All are among the directions represented in $C$ in~\eqref{eq:square-core}, so all edge vectors of the family belong to $\pm C$. In the six-direction case, the same change of basis gives
$$\mathbf p,\quad \mathbf q,\quad \mathbf r,\quad \mathbf p-\mathbf q,\quad \mathbf r+\mathbf p,\quad \mathbf r+\mathbf q,$$
which are all represented in $C$. In both cases, the unused directions of $C$ need not be short.
\end{proof}

For the chain obtained in Proposition~\ref{prop:closed-component}, Lemma~\ref{lem:connected-eight-directions} gives a basis of $\Z^3$ for which all its edge vectors belong to the same set $\pm C$. In these coordinates, Section~\ref{sec:manual-topology} lists the possible tetrahedron translation classes and writes the equations expressing cancellation of their triangular faces. Together with the nonzero determinant sum of this chain, these equations yield a complete family of six co-covering Kuhn tetrahedra for a lattice basis.

\section{Six co-covering Kuhn tetrahedra}\label{sec:manual-topology}

Our goal in this section is to find a basis of $L$ whose six Kuhn tetrahedra are each contained in a translate of $K$. We first list the tetrahedron translation classes with edge vectors in $\pm C$. We then use the cancelling faces and nonzero determinant sum of the chain from Proposition~\ref{prop:closed-component} to obtain the required family.

\subsection{The possible tetrahedron translation classes}

We list all translation classes of four distinct integer points whose pairwise differences lie in $\pm C$. The determinants and reflection pairs of these classes will be used in the next subsection.

Interchanging $\mathbf p$ and $\mathbf q$ preserves $\pm C$, so we may assume
$$\det(\mathbf p,\mathbf q,\mathbf r)=1.$$
We choose new integer coordinates by
$$\mathbf p\mapsto\mathbf e_3,\qquad \mathbf q\mapsto\mathbf e_3-\mathbf e_2,\qquad \mathbf r\mapsto\mathbf e_1-\mathbf e_3.$$
The matrix with the three vectors on the right as columns has determinant one. Thus this coordinate change preserves the integer lattice and the determinant. The eight directions in~\eqref{eq:square-core} are now represented by the following lexicographically positive vectors, listed in increasing order:
\begin{equation}\label{eq:manual-core}
  \begin{gathered}
    \mathbf{v}_1=\mathbf{e}_3,\quad
    \mathbf{v}_2=\mathbf{e}_2-2\mathbf{e}_3,\quad
    \mathbf{v}_3=\mathbf{e}_2-\mathbf{e}_3,\\
    \mathbf{v}_4=\mathbf{e}_2,\quad
    \mathbf{v}_5=\mathbf{e}_1-\mathbf{e}_2,\quad
    \mathbf{v}_6=\mathbf{e}_1-\mathbf{e}_2+\mathbf{e}_3,\\
    \mathbf{v}_7=\mathbf{e}_1-\mathbf{e}_3,\quad
    \mathbf{v}_8=\mathbf{e}_1,\qquad
    C_*=\{\mathbf{v}_1,\ldots,\mathbf{v}_8\}.
  \end{gathered}
\end{equation}
In these coordinates, $\mathbf p=\mathbf v_1$, $\mathbf q=-\mathbf v_3$, and $\mathbf r=\mathbf v_7$. The other five vectors satisfy
$$\mathbf p+\mathbf q=-\mathbf v_2,\quad \mathbf p-\mathbf q=\mathbf v_4,\quad \mathbf r+\mathbf p=\mathbf v_8,\quad \mathbf r+\mathbf q=\mathbf v_5,\quad \mathbf r+\mathbf p+\mathbf q=\mathbf v_6.$$
Hence the coordinate change sends $\pm C$ exactly to $\pm C_*$.

Consider four distinct integer points whose pairwise differences belong to $\pm C_*$. Order the vertices lexicographically and translate the least vertex to zero. Each remaining vertex belongs to $\pm C_*$ and is lexicographically positive, so it belongs to $C_*$. We use the oriented translation classes
$$[ab]=[\mathbf0,\mathbf v_a,\mathbf v_b]\quad(a<b),\qquad[abc]=[\mathbf0,\mathbf v_a,\mathbf v_b,\mathbf v_c]\quad(a<b<c).$$
The increasing indices give the lexicographic order of the vertices. Any reordering multiplies the oriented simplex by the sign of the permutation. Thus every oriented four-point translation class under consideration is represented by some $[abc]$ up to sign.

To determine which triples of indices satisfy the edge conditions, for $a=1,\ldots,8$ define
$$N_a=\{b\in\{1,\ldots,8\}: b>a,\ \mathbf v_b-\mathbf v_a\in\pm C_*\}.$$
Subtracting the vectors in~\eqref{eq:manual-core} gives
\begin{align*}
  N_1&=\{3,4,6,8\},&
  N_2&=\{3,7\},&
  N_3&=\{4,7,8\},&
  N_4&=\{8\},\\
  N_5&=\{6,7,8\},&
  N_6&=\{7,8\},&
  N_7&=\{8\},&
  N_8&=\varnothing.
\end{align*}
For $[abc]$, the three edge vectors from zero are $\mathbf v_a,\mathbf v_b,\mathbf v_c\in C_*$. The other three edge vectors belong to $\pm C_*$ exactly when
$$b,c\in N_a,\qquad c\in N_b.$$
Choosing $b\in N_a$ and then $c\in N_a\cap N_b$ therefore lists all possible classes. The resulting index triples are
$$
  \begin{array}{c|l@{\qquad}c|l}
    a&abc&a&abc\\\hline
    1&134,138,148,168&5&567,568,578\\
    2&237&6&678\\
    3&348,378&7&\varnothing\\
    4&\varnothing&8&\varnothing.
  \end{array}
$$

There are eleven translation classes. In these orientation-preserving coordinates, the normalized determinant from Proposition~\ref{prop:nerve-degree} is
$$d([abc])=\det(\mathbf v_a,\mathbf v_b,\mathbf v_c).$$
It vanishes precisely when the four vertices are coplanar. For a nondegenerate tetrahedron, its sign records the orientation. Computing these determinants shows that only $[134]$ is degenerate. Its vertex set $\{\mathbf0,\mathbf v_1,\mathbf v_3,\mathbf v_4\}$ is a parallelogram because $\mathbf v_1+\mathbf v_3=\mathbf v_4$.

The remaining ten classes are nondegenerate. Reflecting their vertices through the origin and again translating the lexicographically least vertex to zero gives the five pairs in Table~\ref{tab:manual-tetrahedra}. Denote the two oriented representatives in pair $i$ by $T_i^+$ and $T_i^-$, as listed in the table, and put
$$\varepsilon_i=d(T_i^+)=-d(T_i^-), \qquad i=1,\ldots,5.$$

\begin{table}[htbp]
\centering\small
  \caption{The ten nondegenerate translation classes, paired by central reflection.}
\label{tab:manual-tetrahedra}
\begin{tabular}{@{}cccccc@{}}
\toprule
$i$&1&2&3&4&5\\\midrule
$T_i^+/T_i^-$
&$[138]/[678]$
&$[148]/[578]$
&$[168]/[378]$
&$[348]/[568]$
&$[237]/[567]$\\
$\varepsilon_i$&$-1$&$-1$&$1$&$1$&$1$\\
\bottomrule
\end{tabular}
\end{table}

\subsection{Six Kuhn tetrahedra from face cancellation}

We now use face cancellation to obtain the six Kuhn tetrahedra for one basis of $\Z^3$. We state the result for a collection $\mathcal T$ of translation classes of four distinct integer points. Membership in $\mathcal T$ concerns the vertex sets, whereas chains use oriented representatives. The collection is reflection-invariant if reflecting all vertices through the origin preserves membership in $\mathcal T$.

\begin{proposition}\label{prop:manual-six-path-reduction}
Let $\mathcal T$ be reflection-invariant. Suppose every degenerate member is a parallelogram and all pairwise difference vectors of every nondegenerate member belong to $\pm C_*$. If an integral periodic chain $z$ supported on $\mathcal T$ satisfies
$$\partial z=0,\qquad d(z)\ne0,$$
then, for some unimodular matrix $V$, the translation class of every $V\Delta_{ijk}$, $(i,j,k)\in\mathfrak S_3$, belongs to $\mathcal T$.
\end{proposition}

\begin{proof}
Let $\iota$ denote central reflection on formal simplices,
$$\iota[\mathbf x_0,\ldots,\mathbf x_k]=[-\mathbf x_0,\ldots,-\mathbf x_k],$$
and extend it linearly to chains. Reflection commutes with the boundary operator and changes the sign of three-dimensional determinants.

For a parallelogram, translating the reflected vertices by $\mathbf a+\mathbf b$ gives
\begin{align*}
  \iota[\mathbf0,\mathbf a,\mathbf b,\mathbf a+\mathbf b]
  &=[\mathbf a+\mathbf b,\mathbf b,\mathbf a,\mathbf0]\\
  &=[\mathbf0,\mathbf a,\mathbf b,\mathbf a+\mathbf b].
\end{align*}
The last equality uses the even permutation that reverses the order of four vertices. Thus each parallelogram class is unchanged by $\iota$. Define $w=z-\iota z$. By linearity,
\begin{equation}\label{eq:manual-antisymmetrization}
  \partial w=\partial z-\iota(\partial z)=0,
  \quad d(w)=d(z)-d(\iota z)=2d(z)\neq 0.
\end{equation}
The parallelogram coefficients cancel in $w$, so every tetrahedron occurring with nonzero coefficient is nondegenerate. Its vertex-set translation class belongs to $\mathcal T$ because $\mathcal T$ is reflection-invariant.

Central reflection reverses the lexicographic order of four vertices, and restoring that order has sign $(-1)^6=1$. Hence
$$\iota T_i^+=T_i^-,\qquad i=1,\ldots,5.$$
Also $\iota w=-w$, so the coefficients of $T_i^+$ and $T_i^-$ in $w$ are opposite. Put
$$\Xi_i=\varepsilon_i(T_i^+-T_i^-).$$
Since $d(\Xi_i)=2$, we can write
\begin{equation}\label{eq:manual-pair-coordinates}
  w=\sum_{i=1}^5q_i\Xi_i,\qquad
  q_i\in\Z,\qquad d(w)=2\sum_{i=1}^5q_i.
\end{equation}
Whenever $q_i\neq 0$, both $T_i^+$ and $T_i^-$ occur in $w$, so the translation classes of their vertex sets belong to $\mathcal T$.

We next use $\partial w=0$ to relate the coefficients $q_i$. For a class $[abc]$ listed in the preceding subsection, the differences $\mathbf v_b-\mathbf v_a$ and $\mathbf v_c-\mathbf v_a$ belong to $C_*$. Write
$$\mathbf v_b-\mathbf v_a=\mathbf v_r,\qquad \mathbf v_c-\mathbf v_a=\mathbf v_s.$$
Translation preserves lexicographic order, so $r<s$. Translating the face opposite zero by $-\mathbf v_a$ gives
$$\partial[abc]=[rs]-[bc]+[ac]-[ab].$$
For example,
$$\partial[138]=[27]-[38]+[18]-[13].$$

Only three of the face equations are needed. Apply this boundary formula to the five chains $\Xi_i$ and compare the coefficients of the oriented faces $-[27]$, $[48]$, and $-[18]$ in $\partial w$. These three triangular faces represent distinct face classes. Since $\partial w=0$, the total coefficient of each face class must vanish. For $\boldsymbol q=(q_1,\ldots,q_5)^{\mathsf T}$, these equations are
\begin{equation}\label{eq:manual-boundary-matrices}
  \begin{pmatrix}
    1&0&0&0&-1\\
    0&1&0&-1&0\\
    1&1&-1&0&0
  \end{pmatrix}
  \boldsymbol q=0,
  \qquad
  q_5=q_1,\quad q_4=q_2,\quad q_3=q_1+q_2.
\end{equation}
Substituting into~\eqref{eq:manual-pair-coordinates} gives
$$d(w)=2\sum_{i=1}^5q_i=6(q_1+q_2)\neq 0.$$
Thus $q_3\neq 0$, and at least one of $q_1,q_2$ is nonzero. If $q_1\neq 0$, then $q_1,q_3,q_5$ are all nonzero, so both members of pairs $1,3,5$ belong to $\mathcal T$. If $q_2\neq 0$, the same conclusion holds for pairs $2,3,4$.

We identify these two sets of three pairs with the six Kuhn tetrahedra for the bases in Table~\ref{tab:manual-bases}. For $V=[\mathbf a_1\ \mathbf a_2\ \mathbf a_3]$, taking the three basis steps in the order $i,j,k$ gives
$$V\Delta_{ijk} =\{\mathbf0,\mathbf a_i,\mathbf a_i+\mathbf a_j, \mathbf a_1+\mathbf a_2+\mathbf a_3\}.$$

\begin{table}[htbp]
\centering\small
\setlength{\tabcolsep}{5pt}
  \caption{The two bases used to obtain the six Kuhn tetrahedra.}
\label{tab:manual-bases}
\begin{tabular}{@{}ccccc@{}}
\toprule
Basis&$\mathbf a_1$&$\mathbf a_2$&$\mathbf a_3$&$\det V$\\
\midrule
$V_1$&$(1,-1,0)$&$(0,1,-1)$&$(0,0,1)$&$1$\\
$V_2$&$(-1,1,-1)$&$(1,0,-1)$&$(0,0,1)$&$-1$\\
\bottomrule
\end{tabular}
\end{table}

Compute these partial sums and translate each four-point set so that its lexicographically least vertex becomes the origin. The resulting identifications of vertex sets up to translation are
$$
  \begin{array}{c|ccc|c}
    V&V\Delta_{123}&V\Delta_{132}&V\Delta_{213}
    &\text{reflection pairs}\\\hline
    V_1&T_2^-&T_4^-&T_3^-&2,3,4\\
    V_2&T_1^-&T_3^+&T_5^+&1,3,5.
  \end{array}
$$
Here each $T_i^\pm$ denotes its underlying vertex-set class, independently of orientation. For example,
\begin{align*}
  V_2\Delta_{123}=\{\mathbf0,-\mathbf v_6,\mathbf v_2,\mathbf v_3\},\quad V_2\Delta_{123}+\mathbf v_6=\{\mathbf0,\mathbf v_6,\mathbf v_7,\mathbf v_8\},
\end{align*}
which is the vertex set of $T_1^-$.

The three remaining permutations are obtained by reversing the step order. The identity
$$(\mathbf a_1+\mathbf a_2+\mathbf a_3)-V\Delta_{ijk} =V\Delta_{kji}$$
shows that each reverse path is a translate of the central reflection of the original path. Thus the other members of the three reflection pairs give the remaining three Kuhn tetrahedra. Figure~\ref{fig:six-paths} illustrates these pairings in the basis parallelepiped.

Consequently, pairs $2,3,4$ give all six sets $V_1\Delta_{ijk}$, while pairs $1,3,5$ give all six sets $V_2\Delta_{ijk}$. At least one of these two sets of pairs belongs to $\mathcal T$, which proves the proposition.
\end{proof}

\begin{figure}[htbp]
\centering
\begin{tikzpicture}[x=1cm,y=1cm,line cap=round,line join=round,
                    font=\footnotesize,>=stealth]
\begin{scope}[xshift=0.0cm]
  \coordinate (o) at (0,0);
  \coordinate (p1) at (1.55,-0.1);
  \coordinate (p2) at (0.8,0.65);
  \coordinate (p3) at (0,1.7);
  \coordinate (p12) at (2.35,0.55);
  \coordinate (p13) at (1.55,1.6);
  \coordinate (p23) at (0.8,2.35);
  \coordinate (s) at (2.35,2.25);
  \fill[black!8] (o)--(p1)--(p12)--(s)--cycle;
  \fill[black!16] (o)--(s)--(p23)--(p3)--cycle;
  \draw[black!35,line width=.35pt]
    (o)--(p1)--(p12)--(p2)--cycle
    (p3)--(p13)--(s)--(p23)--cycle
    (o)--(p3) (p1)--(p13) (p12)--(s) (p2)--(p23);
  \draw[black!40,line width=.35pt] (o)--(s);
  \draw[->,line width=.85pt] (o)--(p1);
  \draw[->,line width=.85pt] (p1)--(p12);
  \draw[->,line width=.85pt] (p12)--(s);
  \draw[->,densely dashed,line width=.85pt] (o)--(p3);
  \draw[->,densely dashed,line width=.85pt] (p3)--(p23);
  \draw[->,densely dashed,line width=.85pt] (p23)--(s);
  \fill (o) circle (1.2pt);
  \fill (p1) circle (1.2pt);
  \fill (p12) circle (1.2pt);
  \fill (p3) circle (1.2pt);
  \fill (p23) circle (1.2pt);
  \fill (s) circle (1.2pt);
  \node[anchor=north east,inner sep=2pt] at (o) {$\mathbf{0}$};
  \node[anchor=south west,inner sep=2pt] at (s) {$\mathbf{s}$};
  \node[anchor=north,inner sep=2.5pt,fill=white,fill opacity=.85,text opacity=1] at (p1) {$\mathbf{a}_1$};
  \node[anchor=west,inner sep=2.5pt,fill=white,fill opacity=.85,text opacity=1] at (p12) {$\mathbf{a}_1+\mathbf{a}_2$};
  \node[anchor=east,inner sep=2.5pt,fill=white,fill opacity=.85,text opacity=1] at (p3) {$\mathbf{a}_3$};
  \node[anchor=south,inner sep=2.5pt,fill=white,fill opacity=.85,text opacity=1] at (p23) {$\mathbf{a}_2+\mathbf{a}_3$};
  \node at (1.15,-.80) {(a) $123$ (solid), $321$ (dashed)};
\end{scope}
\begin{scope}[xshift=4.55cm]
  \coordinate (o) at (0,0);
  \coordinate (p1) at (1.55,-0.1);
  \coordinate (p2) at (0.8,0.65);
  \coordinate (p3) at (0,1.7);
  \coordinate (p12) at (2.35,0.55);
  \coordinate (p13) at (1.55,1.6);
  \coordinate (p23) at (0.8,2.35);
  \coordinate (s) at (2.35,2.25);
  \fill[black!8] (o)--(p1)--(s)--(p13)--cycle;
  \fill[black!16] (o)--(p2)--(s)--(p23)--cycle;
  \draw[black!35,line width=.35pt]
    (o)--(p1)--(p12)--(p2)--cycle
    (p3)--(p13)--(s)--(p23)--cycle
    (o)--(p3) (p1)--(p13) (p12)--(s) (p2)--(p23);
  \draw[black!40,line width=.35pt] (o)--(s);
  \draw[->,line width=.85pt] (o)--(p1);
  \draw[->,line width=.85pt] (p1)--(p13);
  \draw[->,line width=.85pt] (p13)--(s);
  \draw[->,densely dashed,line width=.85pt] (o)--(p2);
  \draw[->,densely dashed,line width=.85pt] (p2)--(p23);
  \draw[->,densely dashed,line width=.85pt] (p23)--(s);
  \fill (o) circle (1.2pt);
  \fill (p1) circle (1.2pt);
  \fill (p13) circle (1.2pt);
  \fill (p2) circle (1.2pt);
  \fill (p23) circle (1.2pt);
  \fill (s) circle (1.2pt);
  \node[anchor=north east,inner sep=2pt] at (o) {$\mathbf{0}$};
  \node[anchor=south west,inner sep=2pt] at (s) {$\mathbf{s}$};
  \node[anchor=north,inner sep=2.5pt,fill=white,fill opacity=.85,text opacity=1] at (p1) {$\mathbf{a}_1$};
  \node[anchor=north west,inner sep=2.5pt,fill=white,fill opacity=.85,text opacity=1] at (p13) {$\mathbf{a}_1+\mathbf{a}_3$};
  \node[anchor=south east,inner sep=2.5pt,fill=white,fill opacity=.85,text opacity=1] at (p2) {$\mathbf{a}_2$};
  \node[anchor=south,inner sep=2.5pt,fill=white,fill opacity=.85,text opacity=1] at (p23) {$\mathbf{a}_2+\mathbf{a}_3$};
  \node at (1.15,-.80) {(b) $132$ (solid), $231$ (dashed)};
\end{scope}
\begin{scope}[xshift=9.1cm]
  \coordinate (o) at (0,0);
  \coordinate (p1) at (1.55,-0.1);
  \coordinate (p2) at (0.8,0.65);
  \coordinate (p3) at (0,1.7);
  \coordinate (p12) at (2.35,0.55);
  \coordinate (p13) at (1.55,1.6);
  \coordinate (p23) at (0.8,2.35);
  \coordinate (s) at (2.35,2.25);
  \fill[black!8] (o)--(p12)--(s)--cycle;
  \fill[black!16] (o)--(s)--(p3)--cycle;
  \draw[black!35,line width=.35pt]
    (o)--(p1)--(p12)--(p2)--cycle
    (p3)--(p13)--(s)--(p23)--cycle
    (o)--(p3) (p1)--(p13) (p12)--(s) (p2)--(p23);
  \draw[black!40,line width=.35pt] (o)--(s);
  \draw[->,line width=.85pt] (o)--(p2);
  \draw[->,line width=.85pt] (p2)--(p12);
  \draw[->,line width=.85pt] (p12)--(s);
  \draw[->,densely dashed,line width=.85pt] (o)--(p3);
  \draw[->,densely dashed,line width=.85pt] (p3)--(p13);
  \draw[->,densely dashed,line width=.85pt] (p13)--(s);
  \fill (o) circle (1.2pt);
  \fill (p2) circle (1.2pt);
  \fill (p12) circle (1.2pt);
  \fill (p3) circle (1.2pt);
  \fill (p13) circle (1.2pt);
  \fill (s) circle (1.2pt);
  \node[anchor=north east,inner sep=2pt] at (o) {$\mathbf{0}$};
  \node[anchor=south west,inner sep=2pt] at (s) {$\mathbf{s}$};
  \node[anchor=north,inner sep=2.5pt,fill=white,fill opacity=.85,text opacity=1] at (p2) {$\mathbf{a}_2$};
  \node[anchor=west,inner sep=2.5pt,fill=white,fill opacity=.85,text opacity=1] at (p12) {$\mathbf{a}_1+\mathbf{a}_2$};
  \node[anchor=east,inner sep=2.5pt,fill=white,fill opacity=.85,text opacity=1] at (p3) {$\mathbf{a}_3$};
  \node[anchor=north west,inner sep=2.5pt,fill=white,fill opacity=.85,text opacity=1] at (p13) {$\mathbf{a}_1+\mathbf{a}_3$};
  \node at (1.15,-.80) {(c) $213$ (solid), $312$ (dashed)};
\end{scope}
\end{tikzpicture}
  \caption{Six basis-step paths in an affine projection, paired by reflection $\mathbf{x}\mapsto\mathbf{s}-\mathbf{x}$, where $\mathbf{s}=\mathbf{a}_1+\mathbf{a}_2+\mathbf{a}_3$. Before projection, the shaded four-point convex hulls form the Kuhn subdivision of the basis parallelepiped.}
\label{fig:six-paths}
\end{figure}

\subsection{Application to the covering chain}

Let $c$ be the chain supplied by Proposition~\ref{prop:closed-component}. Choose the basis of $\Z^3$ from Lemma~\ref{lem:connected-eight-directions} with positive orientation, and apply the determinant-one coordinate change leading to~\eqref{eq:manual-core}. Continue to denote the resulting lattice basis matrix and integer chain by $B$ and $c$. Every edge vector of $c$ now belongs to $\pm C_*$. These coordinate changes preserve the boundary relation and normalized determinant, so
$$\det B=D>0,\qquad \partial c=0,\qquad d(c)\neq 0.$$

Let $\mathcal T_c$ consist of the translation classes of the vertex sets of tetrahedra occurring in $c$, together with their central reflections. Every member is nondegenerate and has all pairwise difference vectors in $\pm C_*$. Since $K=-K$, central reflection preserves containment in a translate of $K$. Thus every member of $\mathcal T_c$ is co-covering. The collection $\mathcal T_c$ is reflection-invariant and supports $c$, so Proposition~\ref{prop:manual-six-path-reduction} applies.

We obtain a unimodular matrix $V$ such that the translation class of every $V\Delta_{ijk}$, $(i,j,k)\in\mathfrak S_3$, belongs to $\mathcal T_c$. Since these classes are co-covering, their images under $B$ satisfy
\begin{equation}\label{eq:covering-six-paths}
  \forall(i,j,k)\in\mathfrak S_3,\quad
  \exists\,\mathbf t_{ijk}\in\R^3,\qquad
  BV\Delta_{ijk}\subseteq\mathbf t_{ijk}+K.
\end{equation}
Each tetrahedron may have a different containing translate. Since $V$ is unimodular, $BV$ is a basis matrix of $L$ and
$$|\det(BV)|=D.$$
Theorem~\ref{thm:packet-radius} in the next section will use these six containments to give $D\leq 4$.

\section{A sharp determinant bound for six Kuhn tetrahedra}\label{sec:six-path}

We prove that if each of the six sets $B\Delta_{ijk}$ is contained in a translate of $rO$, then
$$|\det B|\le\frac{32}{27}r^3.$$
Applied to the basis matrix $BV$ in~\eqref{eq:covering-six-paths} with $r=3/2$, this bound will give $D\le4$.

For a real $3\times3$ matrix $B$, define
\begin{equation}\label{eq:packet-radius-definition}
  \rho_{\mathrm p}(B)
  =\max_{(i,j,k)\in\mathfrak S_3}
  \min_{\mathbf z\in\R^3}
  \max_{\mathbf v\in\Delta_{ijk}}
  |B\mathbf v-\mathbf z|_1.
\end{equation}
Thus $\rho_{\mathrm p}(B)$ is the least radius $r$ such that each of the six sets $B\Delta_{ijk}$ is contained in a translate of $rO$. The containing translates may have different centers. For each set, the function being minimized is continuous and tends to infinity as $|\mathbf z|_1\to\infty$, so the minimum is attained.

Put
\begin{equation}\label{eq:sp-change}
  T=
  \begin{pmatrix}
    1&1&-1\\
    1&-1&1\\
    -1&1&1
  \end{pmatrix},
  \qquad
  T^{-1}=\frac12
  \begin{pmatrix}
    1&1&0\\
    1&0&1\\
    0&1&1
  \end{pmatrix}.
\end{equation}
Set $A=BT^{-1}$. The columns of $A$ are the half-sums of the three pairs of columns of $B$. This change puts the vertices of the Kuhn tetrahedra into the paired form used in Section~\ref{subsec:six-containments}. Since $B=AT$ and $|\det T|=4$,
$$|\det B|=4|\det A|.$$

Let $\sigma_1(A),\sigma_2(A),\sigma_3(A)$ be the singular values of $A$, namely the nonnegative square roots of the eigenvalues of $A^{\mathsf T}A$. Their sum is the {\it nuclear norm}
$$\|A\|_*=\sum_{i=1}^3\sigma_i(A).$$
The arithmetic--geometric mean inequality gives
$$|\det B|=4\prod_{i=1}^3\sigma_i(A) \leq 4\left(\frac{\|A\|_*}{3}\right)^3.$$
It therefore suffices to prove $\|A\|_*\le2\rho_{\mathrm p}(B)$. The following theorem states this estimate and the resulting determinant bound.

\begin{theorem}\label{thm:packet-radius}
For every real $3\times3$ matrix $B$,
\begin{equation}\label{eq:packet-nuclear}
  \rho_{\mathrm p}(B)\geq \frac12\|BT^{-1}\|_*,
  \qquad
  \rho_{\mathrm p}(B)^3\geq \frac{27}{32}|\det B|.
\end{equation}
Both constants are best possible, with equality at $B=T$.
\end{theorem}

We first derive bounds on the sum of the absolute entries in each $2\times2$ block of $A$ and the entry in the remaining row and column. We then prove a matrix inequality that uses these bounds to control $\|A\|_*$.

\subsection{Inequalities from the six containments} \label{subsec:six-containments}

Each containing translate gives four $\ell_1$ distance bounds, one for each vertex. We combine these bounds to obtain inequalities for the entries of $A=BT^{-1}$ that do not involve the center.

\begin{lemma}\label{lem:six-containment-five-term}
Let $B\in\R^{3\times3}$ and $r\ge0$, and put $A=BT^{-1}=(a_{ij})$. Suppose that each set $B\Delta_{ijk}$, $(i,j,k)\in\mathfrak S_3$, is contained in a translate of $rO$. For any two-element subsets $I,J\subset\{1,2,3\}$, write
$$\{i_0\}=\{1,2,3\}\setminus I,\qquad \{j_0\}=\{1,2,3\}\setminus J.$$
Then
$$\sum_{\substack{i\in I\\j\in J}}|a_{ij}|+|a_{i_0j_0}|\le2r.$$
\end{lemma}

\begin{proof}
We first show that row and column permutations of $A$ preserve the six containment conditions. This will reduce the proof to $I=J=\{1,2\}$.

A row permutation of $A$ induces the same coordinate permutation of $B=AT$. Since coordinate permutations preserve $O$, the containment conditions remain valid.

The columns of $T$ are precisely the three vectors with two entries equal to $1$ and one entry equal to $-1$. A row permutation of $T$ therefore permutes its columns. Consequently, for every permutation matrix $R$, there is a permutation matrix $Q$ such that
$$RT=TQ,\qquad (AR)T=BQ.$$
If $Q\mathbf e_i=\mathbf e_{\pi(i)}$, then
$$BQ\Delta_{ijk}=B\Delta_{\pi(i)\pi(j)\pi(k)}.$$
Thus a column permutation of $A$ merely permutes the six sets under consideration and also preserves the hypothesis.

We now consider $I=J=\{1,2\}$. Direct calculation gives
$$T\Delta_{312}-\mathbf e_2=\{\pm\mathbf e_2,\ \mathbf e_3\pm\mathbf e_1\}.$$
Choose $\mathbf c$ such that $B\Delta_{312}\subseteq\mathbf c+rO$, and put
$$\mathbf z=\mathbf c-A\mathbf e_2=(z_1,z_2,z_3).$$
Since $AT=B$, we have
\begin{align*}
  A\{\pm\mathbf e_2,\ \mathbf e_3\pm\mathbf e_1\}
  &=B\Delta_{312}-A\mathbf e_2\\
  &\subseteq\mathbf z+rO.
\end{align*}
Adding the four resulting distance bounds and using
$$|u+v|+|u-v|=2\max\{|u|,|v|\}$$
gives
\begin{align*}
  4r &\geq \sum_{\epsilon=\pm1}\bigl(|\epsilon A\mathbf e_2-\mathbf z|_1 +|A(\mathbf e_3+\epsilon\mathbf e_1)-\mathbf z|_1\bigr)\\
  &=2\sum_{i=1}^3\bigl(\max\{|a_{i2}|,|z_i|\}+\max\{|a_{i1}|,|a_{i3}-z_i|\}\bigr)\\
  &\geq2\bigl(|a_{11}|+|a_{12}|+|a_{21}|+|a_{22}|+|z_3|+|a_{33}-z_3|\bigr)\\
  &\geq 2\bigl(|a_{11}|+|a_{12}|+|a_{21}|+|a_{22}|+|a_{33}|\bigr).
\end{align*}
For $i=1,2$, we bound the two maxima below by $|a_{i2}|$ and $|a_{i1}|$. For $i=3$, we instead use $|z_3|$ and $|a_{33}-z_3|$. The last inequality follows from $|z_3|+|a_{33}-z_3|\geq|a_{33}|$.

Dividing by two proves the required bound for $I=J=\{1,2\}$. The row and column permutations considered above give the result for every choice of $I$ and $J$.
\end{proof}

\subsection{A bound for the nuclear norm}

We use the bounds in Lemma~\ref{lem:six-containment-five-term} to control $\|BT^{-1}\|_*$. For this purpose, define
$$
  \mathcal F_5=\left\{R_1
  \begin{pmatrix}
    a&b&0\\
    c&-abc&0\\
    0&0&d
  \end{pmatrix}
  R_2^{\mathsf T}
  \ \middle|\
  \begin{array}{l}
    a,b,c,d\in\{\pm1\},\\
    R_1,R_2\text{ are }3\times3\text{ permutation matrices}
  \end{array}
  \right\}.
$$
Each member has four entries in $\{\pm1\}$ in a $2\times2$ block, with product $-1$, and one further entry in $\{\pm1\}$ in the remaining row and column. All other entries are zero. Thus
$$|\mathcal F_5|=9\cdot8\cdot2=144.$$
For real $3\times3$ matrices $X=(x_{ij})$ and $Y=(y_{ij})$, write
$$\langle X,Y\rangle=\sum_{i,j=1}^3x_{ij}y_{ij}=\operatorname{tr}(X^{\mathsf T}Y),$$
and define
$$h_5(A)=\max_{F\in\mathcal F_5}\langle F,A\rangle.$$
Thus $h_5(A)$ is the largest of these $144$ signed five-term sums. The set $\mathcal F_5$ is centrally symmetric. It is also invariant under transposition and independent signed row and column permutations, since these operations preserve the block structure and the product of the four block entries.

Under the hypotheses of Lemma~\ref{lem:six-containment-five-term}, put $A=BT^{-1}=(a_{ij})$. For $F\in\mathcal F_5$, let $I,J$ be the row and column sets of its $2\times2$ block, and let $i_0,j_0$ be the remaining row and column indices. The lemma gives
$$\langle F,A\rangle\leq \sum_{\substack{i\in I\\j\in J}}|a_{ij}|+|a_{i_0j_0}|\leq 2r.$$
Taking the maximum over $F\in\mathcal F_5$, we obtain
$$h_5(BT^{-1})\leq2r.$$
It remains to prove
$$\|A\|_*\le h_5(A)$$
for every real $3\times3$ matrix $A$.

The singular-value decomposition gives the variational formula
$$\|A\|_*=\max_{R\in\mathrm O(3)}\langle R,A\rangle, \qquad\mathrm O(3)=\{R\in\R^{3\times3}:R^{\mathsf T}R=I\}.$$
It therefore suffices to show that
$$\mathrm O(3)\subseteq\conv\mathcal F_5.$$
We establish this containment using the known descriptions of the Clifford and Mermin polytopes.

Let
$$\mathcal G=\{C\in\{0,\pm1\}^{3\times3}: C^{\mathsf T}C=I,\ \det C=1\},\qquad|\mathcal G|=24.$$
These are the signed permutation matrices of determinant one. Their convex hull is the Clifford polytope studied in~\cite{BCLLSU2006}. Consider its polar
$$P=(\conv\mathcal G)^\circ=\{X\in\R^{3\times3}:\langle C,X\rangle\leq 1 \text{ for every }C\in\mathcal G\},$$
and put
$$\mathcal B=\{F\in\mathcal F_5:\det F=-2\}.$$
Two vertices of $P$ are called {\it neighbors} if an edge of $P$ joins them. The following result describes the vertices of $P$ and the neighbors of those outside $\mathcal B$.

\begin{lemma}[Buhrman et al.~\cite{BCLLSU2006} and Okay--Chung--Ipek~\cite{OCI2023}]\label{lem:clifford-mermin-description}
The polytope $P$ has $120$ vertices. They consist of the $72$ matrices in $\mathcal B$ and the $48$ matrices whose nonzero entries are three arbitrary signs in one full row or column. If $E$ is a vertex of the latter type, every neighbor $F$ belongs to $\mathcal B$ and satisfies
$$\langle E,F\rangle=2.$$
\end{lemma}

We recall how the cited results give this description in our coordinates. The vertex list follows by polarity from the Clifford facet description in~\cite[Sec.~4, Eq.~(8)]{BCLLSU2006}. In the coefficient-matrix coordinates of~\cite[Lemma~2]{OCI2023}, the polytope $P$ is the Mermin polytope $\mathrm{MP}_1$. Its vertex classification in~\cite[Theorem~3]{OCI2023} identifies the $48$ three-entry vertices as type~1 and the $72$ members of $\mathcal B$ as type~2.

By~\cite[Theorem~5 and Lemma~15]{OCI2023}, every type~1 vertex has only type~2 neighbors, and the adjacent pairs of these two types form one symmetry orbit. The adjacent representatives in~\cite[Eqs.~(34)--(35)]{OCI2023} have coefficient matrices
$$
  E_0=
  \begin{pmatrix}
    0&1&0\\
    0&-1&0\\
    0&1&0
  \end{pmatrix},
  \qquad
  F_0=
  \begin{pmatrix}
    1&1&0\\
    1&-1&0\\
    0&0&1
  \end{pmatrix},
  \qquad
  \langle E_0,F_0\rangle=2.
$$
The symmetries used there act by signed row and column permutations and transposition. These operations preserve the inner product, so every such adjacent pair satisfies $\langle E,F\rangle=2$.

We now remove the $48$ three-entry vertices by intersecting $P$ with half-spaces whose boundary hyperplanes pass through their neighbors. The resulting description of $\conv\mathcal B$ will give the required containment of orthogonal matrices.

\begin{lemma}\label{lem:five-term-nuclear}
Every real $3\times3$ matrix $A$ satisfies
$$\|A\|_*\leq h_5(A).$$
\end{lemma}

\begin{proof}
Fix a vertex $E$ of $P$ with exactly three nonzero entries, and intersect $P$ with the half-space
$$\{X\in\R^{3\times3}:\langle E,X\rangle\leq2\}.$$
Since $\langle E,E\rangle=3$, this intersection excludes $E$. Every $F\in\mathcal B$ has at most two nonzero entries in the row or column occupied by $E$, each of absolute value one. Hence
$$\langle E,F\rangle\leq 2,$$
so all vertices in $\mathcal B$ are retained.

Let $E'$ be another three-entry vertex. If its nonzero entries occupy the same positions as those of $E$, at least one sign differs. Otherwise, the two sets of positions overlap in at most one entry. In both cases,
$$\langle E,E'\rangle\leq 1<2.$$
Thus $E$ is the only excluded vertex.

By Lemma~\ref{lem:clifford-mermin-description}, every neighbor of $E$ lies on the hyperplane $\langle E,X\rangle=2$. A new vertex could arise only where this hyperplane meets the interior of an edge joining $E$ to a retained vertex. Every such edge instead meets the hyperplane at its existing endpoint. The intersection therefore creates no new vertices. Figure~\ref{fig:vertex-truncation} illustrates this argument in two dimensions.

Each remaining three-entry vertex satisfies the added inequality strictly, and all its neighbors are retained. Its incident edges are therefore unchanged. We may repeat the same argument at each of the remaining three-entry vertices. After all $48$ intersections, the remaining vertices are exactly the members of $\mathcal B$.

\begin{figure}[htbp]
\centering
\begin{tikzpicture}[
  x=1.15cm,y=1.15cm,
  line cap=round,line join=round,
  font=\small
]
  \coordinate (E) at (0,1.73205);
  \coordinate (F1) at (-1.55,1.15470);
  \coordinate (F2) at (1.55,1.15470);
  \coordinate (L) at (-2.05,.1);
  \coordinate (BL) at (-1.2,-.65);
  \coordinate (BR) at (1.2,-.65);
  \coordinate (R) at (2.05,.1);

  \filldraw[
    fill=black!7,draw=black!70,line width=.6pt
  ]
    (F1)--(L)--(BL)--(BR)--(R)--(F2)--cycle;

  \draw[densely dashed,black!65,line width=.6pt]
    (F1)--(E)--(F2);

  \draw[densely dashed,black!65,line width=.45pt]
    (-2.75,1.15470)--(F1)
    (F2)--(2.75,1.15470);

  \draw[line width=.85pt] (F1)--(F2);

  \foreach \point in {E,F1,F2,L,BL,BR,R}
    \fill (\point) circle (1.5pt);

  \node[anchor=south,inner sep=3pt] at (E)
    {$E,\quad\langle E,E\rangle=3$};
  \node[anchor=south east,inner sep=3pt] at (F1)
    {$F_1$};
  \node[anchor=south west,inner sep=3pt] at (F2)
    {$F_2$};
  \node[anchor=west,inner sep=3pt] at (2.75,1.15470)
    {$\langle E,X\rangle=2$};
  \node at (0,.2)
    {$\langle E,X\rangle\le2$};
\end{tikzpicture}
  \caption{A two-dimensional analogue of the nine-dimensional vertex cut, not a face or projection of $P$. The cutting line passes through the existing neighbors $F_1,F_2$. Removing $E$ introduces no new vertices.}
\label{fig:vertex-truncation}
\end{figure}

For a fixed row $i$, the eight inequalities corresponding to its sign choices are equivalent to
$$\max_{\varepsilon_1,\varepsilon_2,\varepsilon_3\in\{\pm1\}}\sum_{j=1}^3\varepsilon_jx_{ij} =\sum_{j=1}^3|x_{ij}|\leq 2.$$
The same argument applies to each column. Combining these inequalities with the defining inequalities of $P$, we obtain
\begin{equation}\label{eq:five-term-polytope}
  \conv\mathcal B=\left\{X=(x_{ij})\in\R^{3\times3}:
  \begin{array}{ll}
    \langle C,X\rangle\leq 1
    &(C\in\mathcal G),\\
    \sum_{j=1}^3|x_{ij}|\leq 2
    &(i=1,2,3),\\
    \sum_{i=1}^3|x_{ij}|\leq 2
    &(j=1,2,3)
  \end{array}
  \right\}.
\end{equation}

We next show that $-R$ satisfies these inequalities for every three-dimensional rotation $R$. Write
$$\mathrm{SO}(3)=\{R\in\mathrm O(3):\det R=1\},$$
and take $R=(r_{ij})\in\mathrm{SO}(3)$. For $C\in\mathcal G$, the matrix $C^{\mathsf T}R$ is also a rotation. If its rotation angle is $\theta$, its trace is $1+2\cos\theta$. Consequently,
$$\langle C,-R\rangle=-\operatorname{tr}(C^{\mathsf T}R)=-(1+2\cos\theta)\leq 1.$$
Every row and column of $R$ has Euclidean norm one. The Cauchy--Schwarz inequality gives
$$\sum_{j=1}^3|r_{ij}|\le\sqrt3<2 \qquad(i=1,2,3),$$
and
$$\sum_{i=1}^3|r_{ij}|\leq \sqrt3<2 \qquad(j=1,2,3).$$
Thus $-R$ satisfies all the inequalities in~\eqref{eq:five-term-polytope}, proving
$$-\mathrm{SO}(3) \subseteq\conv\mathcal B \subseteq\conv\mathcal F_5.$$
Since $\mathcal F_5$ is centrally symmetric, $\mathrm{SO}(3)$ also lies in $\conv\mathcal F_5$. Together with $\mathrm O(3)=\mathrm{SO}(3)\cup(-\mathrm{SO}(3)),$ this gives
$$\mathrm O(3)\subseteq\conv\mathcal F_5.$$
The variational formula for the nuclear norm and linearity of the inner product now yield
\begin{align*}
  \|A\|_* &=\max_{R\in\mathrm O(3)}\langle R,A\rangle\\
  &\le\max_{X\in\conv\mathcal F_5}\langle X,A\rangle\\
  &=\max_{F\in\mathcal F_5}\langle F,A\rangle=h_5(A).
\end{align*}
\end{proof}

Numerical evidence for an equivalent formulation of $-\mathrm{SO}(3)\subseteq\conv\mathcal B$ was reported in~\cite[Example~2 and Table~III]{RBTOZ2020}.

\subsection{The determinant bound and sharpness}

\begin{proof}[Proof of Theorem~\ref{thm:packet-radius}]
Put $A=BT^{-1}$. Since the minima in~\eqref{eq:packet-radius-definition} are attained, each $B\Delta_{ijk}$ is contained in a translate of $\rho_{\mathrm p}(B)O$. Applying Lemma~\ref{lem:six-containment-five-term} to these six containments, and then Lemma~\ref{lem:five-term-nuclear}, we obtain
$$\|A\|_*\leq h_5(A)\leq 2\rho_{\mathrm p}(B).$$
This proves the first inequality in~\eqref{eq:packet-nuclear}.

Let $\sigma_1,\sigma_2,\sigma_3$ be the singular values of $A$. Since $B=AT$ and $|\det T|=4$, the arithmetic--geometric mean inequality gives
\begin{align*}
  |\det B|&=4|\det A|=4\sigma_1\sigma_2\sigma_3\\
  &\leq 4\left(\frac{\sigma_1+\sigma_2+\sigma_3}{3}\right)^3\\
  &= 4\left(\frac{\|A\|_*}{3}\right)^3
  \leq \frac{32}{27}\rho_{\mathrm p}(B)^3.
\end{align*}
This proves the second inequality.

To verify sharpness, take $B=T$. For the Kuhn tetrahedron $\Delta_{123}$, we have
$$T\Delta_{123}=\{\mathbf0,(1,1,-1),(2,0,0),(1,1,1)\}.$$
With $\mathbf c=(1,\tfrac12,0)$, direct calculation gives
$$|\mathbf x-\mathbf c|_1=\frac{3}{2} \qquad (\mathbf x\in T\Delta_{123}).$$
Thus $T\Delta_{123}\subseteq\mathbf c+\frac{3}{2}O$.

For any $(i,j,k)\in\mathfrak S_3$, let $Q$ be the permutation matrix with columns $\mathbf e_i,\mathbf e_j,\mathbf e_k$. Then
$$\Delta_{ijk}=Q\Delta_{123}.$$
The permutation property of $T$ used in Section~\ref{subsec:six-containments} gives a permutation matrix $P$ such that $TQ=PT$. Since $PO=O$, it follows that
$$T\Delta_{ijk}=TQ\Delta_{123}=PT\Delta_{123} \subseteq P\mathbf c+\frac{3}{2}O.$$
All six Kuhn tetrahedra therefore have containing translates of $\frac{3}{2}O$, so
$$\rho_{\mathrm p}(T)\leq\frac{3}{2}.$$
The first inequality already proved gives the reverse bound
$$\rho_{\mathrm p}(T)\geq \frac{1}{2}\|TT^{-1}\|_*=\frac{1}{2}\|I\|_*=\frac{3}{2}.$$
Hence $\rho_{\mathrm p}(T)=3/2$, with equality in the first inequality of~\eqref{eq:packet-nuclear}. Since $|\det T|=4$, we also have
$$\rho_{\mathrm p}(T)^3=\frac{27}{8}=\frac{27}{32}|\det T|.$$
Thus equality holds in both inequalities at $B=T$, and both constants are best possible.
\end{proof}

Each column of $T$ belongs to the lattice $L_0$ defined in~\eqref{eq:bcc-lattice}, so $T\Z^3\subseteq L_0$. Both lattices have determinant $4$, and consequently
$$T\Z^3=L_0.$$

\section{Completion of the proof}\label{sec:verification}

We combine the six containments in~\eqref{eq:covering-six-paths} with Theorem~\ref{thm:packet-radius} to complete the proof of Theorem~\ref{thm:main}. We then derive an asymptotic consequence for Abelian Cayley graphs with three generators.

\begin{proof}[Proof of Theorem~\ref{thm:main}]
For $K=\frac{3}{2}O$, the lattice $L_0$ in~\eqref{eq:bcc-lattice} satisfies
$$K+L_0=\R^3,\qquad\det L_0=4.$$
Since lattice covering density is invariant under scaling,
$$\theta^l(O)=\theta^l(K)\leq \frac{\vol(K)}{\det L_0}=\frac{9}{8}.$$

For the reverse inequality, it suffices by~\eqref{eq:determinant-goal} to prove that every lattice $L$ with $K+L=\R^3$ satisfies $\det L\leq 4$. Suppose, for a contradiction, that
$$K+L=\R^3, \qquad D=\det L>4.$$
By~\eqref{eq:covering-six-paths}, $L$ has a basis matrix $\widetilde B$ such that each of the six Kuhn tetrahedra $\widetilde B\Delta_{ijk}$ is contained in a translate of $K=\frac{3}{2}O$. Therefore,
$$|\det\widetilde B|=D, \qquad \rho_{\mathrm p}(\widetilde B)\leq \frac{3}{2}.$$
Theorem~\ref{thm:packet-radius} now gives
\begin{align*}
  D=|\det\widetilde B| &\leq\frac{32}{27}\rho_{\mathrm p}(\widetilde B)^3\\
  &\leq  \frac{32}{27}\left(\frac32\right)^3=4,
\end{align*}
contrary to $D>4$. Thus every lattice covering by $K$ has density at least $\vol(K)/4=9/8$. Together with the covering by $L_0$, this proves $\theta^l(O)=9/8$.
\end{proof}

The relation between lattice coverings and Cayley graph diameters established by Dougherty and Faber~\cite{DF2004} gives an application of Theorem~\ref{thm:main}. The degree--diameter problem asks for the largest number of vertices in a graph with prescribed upper bounds on its degree and diameter. Here we consider undirected Abelian Cayley graphs with three generators and their inverses.

Let $H$ be a finite Abelian group generated by $g_1,g_2,g_3$. The Cayley graph
$$G=\operatorname{Cay}\bigl(H,\{\pm g_1,\pm g_2,\pm g_3\}\bigr)$$
has vertex set $H$, with distinct vertices $h,h'$ adjacent when $h'-h=\pm g_i$ for some $i$. Its degree is at most six. Its diameter is the maximum graph distance between two vertices.

\begin{corollary}\label{cor:graph-asymptotic}
Let $M_3(k)$ be the largest order of a finite Abelian group admitting an undirected Cayley graph of diameter at most $k$ with three generators and their inverses. Then
\begin{equation}\label{eq:graph-asymptotic}
  M_3(k)=\frac{32}{27}k^3+O(k^2)
  \qquad(k\to\infty).
\end{equation}
\end{corollary}

\begin{proof}
For the graph $G$ defined above, consider the surjective homomorphism
$$
  \begin{aligned}
    \pi:\Z^3&\longrightarrow H,\\
    (m_1,m_2,m_3)&\longmapsto m_1g_1+m_2g_2+m_3g_3,
  \end{aligned}
  \qquad
  \Lambda=\ker\pi.
$$
Then
$$H\cong\Z^3/\Lambda, \qquad |H|=[\Z^3:\Lambda].$$
Each edge of $G$ adds or subtracts one generator. A vector $\mathbf m=(m_1,m_2,m_3)$ with $\pi(\mathbf m)=h$ gives a walk from $0$ to $h$ of length at most $|\mathbf m|_1$. Conversely, any walk from $0$ to $h$ gives such a vector whose $\ell_1$ norm is at most the length of the walk. Hence
$$\operatorname{dist}_G(0,h)=\min\{|\mathbf m|_1:\mathbf m\in\Z^3,\ \pi(\mathbf m)=h\}.$$

Translations in $H$ preserve graph distance, so
$$\operatorname{diam}(G)=\max_{h\in H}\operatorname{dist}_G(0,h).$$
Consequently, for every integer $k\geq 0$,
 \begin{align*}
    \operatorname{diam}(G)\leq k &\quad\Longleftrightarrow\quad \pi(kO\cap\Z^3)=H\\
    &\quad\Longleftrightarrow\quad (kO\cap\Z^3)+\Lambda=\Z^3.
  \end{align*}
The last condition means that each coset of $\Lambda$ has a representative in $kO\cap\Z^3$.

Conversely, every finite-index sublattice $\Lambda\subseteq\Z^3$ gives such a Cayley graph on $\Z^3/\Lambda$, using the classes of $\mathbf e_1,\mathbf e_2,\mathbf e_3$ and their inverses. Thus $M_3(k)$ is the largest index $[\Z^3:\Lambda]$ among the integer lattices for which $(kO\cap\Z^3)+\Lambda=\Z^3$.

By Theorem~\ref{thm:main}, the optimal lattice covering efficiency of $O$ is
$$\varepsilon=\frac{1}{\theta^l(O)}=\frac{8}{9}.$$
Dougherty and Faber's asymptotic relation in~\cite[Corollary~4.7(a)]{DF2004} gives
$$M_3(k)=\vol(O)\,\varepsilon k^3+O(k^2).$$
Since $\vol(O)=4/3$, we conclude that
$$M_3(k)=\frac43\cdot\frac89\,k^3+O(k^2)=\frac{32}{27}k^3+O(k^2).$$
\end{proof}

\end{document}